\documentclass[11pt,twoside]{article}
\usepackage[margin=1.3in]{geometry}

\let\oldbibliography\thebibliography
\renewcommand{\thebibliography}[1]{%
\oldbibliography{#1}%
\setlength{\itemsep}{0pt}%
}

\usepackage[T1]{fontenc}
\usepackage[utf8]{inputenc}
\usepackage{amsmath,amssymb,amsthm,mathtools}
\usepackage{newtxtext,newtxmath,microtype}
\usepackage[hidelinks]{hyperref}
\allowdisplaybreaks[2]
\numberwithin{equation}{section}
\newtheorem{theorem}{Theorem}[section]
\newtheorem{lemma}[theorem]{Lemma}
\newtheorem{proposition}[theorem]{Proposition}
\newtheorem{corollary}[theorem]{Corollary}
\theoremstyle{remark}\newtheorem{remark}[theorem]{Remark}
\newcommand{\R}{\mathbb R}\newcommand{\Z}{\mathbb Z}
\newcommand{\dd}{\,\mathrm d}
\DeclareMathOperator{\supp}{supp}

\DeclareMathOperator{\Lip}{Lip}

\numberwithin{equation}{section}  

\begin{document} 
\title{ {\bf Failure of superharmonicity near an isolated singularity of the critical biharmonic equation} }

\author{Juncheng Wei \qquad Hui Yang \qquad Ronghao Yang}
\date{\today}
\maketitle

\begin{center}
\begin{minipage}{130mm}
\begin{center}{\bf Abstract}\end{center}
Let $n\ge5$ and $u$ be  a positive singular  solution of
$\Delta^2u=u^{(n+4)/(n-4)}$ in a punctured ball. A long-standing question is where or not $ -\Delta u>0$ in the punctured ball. In this paper we show that the answer is negative. More precisely we construct a positive singular solution whose singularity at the
origin is nonremovable, while its lower limit there is zero. On a normal
subspace, $-\Delta u$ converges to a strictly negative constant. The
spherical minimum has order $r^\beta$, where $2<\beta<4$, and the
spherical Harnack ratio is unbounded. The proof uses
balanced planar polygons, finitely many exterior parameters and a
countable gluing construction. 

\vskip0.10in


\end{minipage}

\end{center}

\section{Introduction}\label{intro:section}
Let $B_r=\{X\in\R^n:|X|<r\}$. We consider
\begin{equation}\label{intro:pde}
\Delta^2u=u^{\frac{n+4}{n-4}}\qquad\text{in }B_1\setminus\{0\},
\end{equation}
where $n\ge5$ is an integer. In \cite{JX}, Jin and Xiong assumed that the positive solution $u\in C^{4}(B_1 \setminus\{0\})$ of \eqref{intro:pde} satisfies 
\begin{equation}\label{eq:all-posi01}
-\Delta u(X)\geq 0
\quad \text{in }B_1 \setminus\{0\}. 
\end{equation} 
Then they proved that either the origin is a removable singularity of $u$, or there exists $C>0$ such that 
\[
\frac{1}{C}|X|^{-\frac{n-4}{2}} \leq u(X) \leq C |X|^{-\frac{n-4}{2}} \quad \text{for } 0< |X|< \frac{1}{2}.
\]
This extends the celebrated result of Caffarelli-Gidas-Spruck \cite{CGS} for the second-order Yamabe equation. Du and Yang~\cite{DY} studied
corresponding local estimates for more general singular sets, also under
lower-order sign assumptions. It is natural to ask: does the sign condition \eqref{eq:all-posi01} hold automatically for singular solutions? In other words, if the origin is a nonremovable singularity, i.e., $\limsup_{X\to 0} u(X) = +\infty$, 
whether there always exists $\rho\in(0,1)$ such that
\begin{equation}\label{eq:all-posi02}
-\Delta u(X)\geq 0
\quad \text{for }0<|X|<\rho. 
\end{equation}
For the equation \eqref{intro:pde} on the entire space $\mathbb{R}^n$, Wei and Xu \cite{WX} have proved that any positive solution always satisfies this sign condition automatically (also known as the superharmonic property). For global singular solutions $u\in C^{4} (\mathbb{R}^n\setminus\{0\})$, similar results can be found in \cite{AGHW,DY,Lin}. However, these proofs rely on the fact that solutions are defined on $\mathbb{R}^n$ or $\mathbb{R}^n\setminus\{0\}$, and hence cannot be applied to the local solutions considered here. The classification of singular positive solutions on $\R^n\setminus\{0\}$ was established by Frank and K\"onig~\cite{FK} based on the works of Lin \cite{Lin} and Guo-Huang-Wang-Wei \cite{GHWW}. 
For the subcritical biharmonic Lane-Emden equation, Huang-Li-Yang \cite{HLY} obtained the superharmonic property for local singular solutions. Nevertheless, their proof is not applicable to the Sobolev critical case.

In this paper, we show that positivity and nonremovability do not imply the sign assumption for the critical biharmonic equation \eqref{intro:pde}, even in a smaller punctured neighborhood. Write $X=(x,z)\in\R^2\times\R^{n-2}$.

\begin{theorem}\label{intro:main}
For every integer $n\ge5$, there exist
$u\in C^\infty(B_1\setminus\{0\})$ and $c_*>0$ such that
$u>0$, $u$ solves \eqref{intro:pde}, and
\begin{align}
\liminf_{X\to0}u(X)&=0,&
\limsup_{X\to0}u(X)&=+\infty,\label{intro:oscillation}\\
\lim_{z\to0,\ z\ne0}(-\Delta u)(0,z)&=-c_*.
\label{intro:negative}
\end{align}
In particular, the origin is the only nonremovable singularity of $u$,
and $u$ is not superharmonic in any punctured neighborhood of the origin.
\end{theorem}

For $0<r<1$, define
\[
m_u(r)=\min_{|X|=r}u(X),\qquad M_u(r)=\max_{|X|=r}u(X).
\]
\begin{corollary}\label{intro:envelope}
The solution in Theorem~\ref{intro:main} can be chosen so that
\begin{equation}\label{intro:minimum}
cr^\beta\le m_u(r)\le Cr^\beta\qquad(0<r<r_*),
\end{equation}
where $c,C,r_*>0$, with $3<\beta<4$ if $n=5$ and
$2<\beta<3$ if $n\ge6$. There is a sequence $r_i\downarrow0$
such that
\begin{equation}\label{intro:maximum}
M_u(r_i)\ge c r_i^{-(n-4+\beta)},\qquad
\frac{M_u(r_i)}{m_u(r_i)}\ge c r_i^{-(n-4+2\beta)}.
\end{equation}
Consequently, for every $\gamma<\beta$, no positive constants $b,r_b$
satisfy $u(X)\ge b|X|^\gamma$ for $0<|X|<r_b$.
The upper estimate $u(X)\le C|X|^{-(n-4)/2}$ also fails near the origin.
\end{corollary}

We describe the construction in Theorem \ref{intro:main} in more detail. Put $s=n-4$, $a=s/2$
and $p=1+8/s$. In every dimension, the sources lie on regular
polygons in the same plane $\{z=0\}$. Let $N$ be even, let $v_\nu$
be the vertices of the unit regular $N$-gon, and set
\[
 R_i=r_0t^i,\qquad r_0=\sqrt t,\qquad
 q_{i\nu}=R_i(v_\nu,0),\qquad M_i=R_i^{s+\beta},\qquad i\in\Z.
\]
Here $0<t<1$ and $2<\beta<4$ will be determined by balancing the
sources. With $T_2$ denoting the Taylor polynomial of degree two at
the origin, the point potential is
\[
 F(X)=\sum_{i\in\Z}\frac{M_i}{N}
       \sum_{\nu=0}^{N-1}(I-T_2)|X-q_{i\nu}|^{-s}.
\]
The subtracted polynomials are biharmonic and therefore leave all
point masses positive. Opposite vertices cancel odd Taylor terms,
so the outer-shell remainder starts at degree four. The restrictions
$\beta<4$ and $\beta>2$ ensure, respectively, summability of this
remainder and of the inner-shell quadratic terms.
The bilateral array is useful because shifting the shell index
amounts to a spatial dilation. Consequently, one balanced polygon
determines the regular jets on every shell, including those that
accumulate at the origin. The exterior shells provide the freedom
needed for this exact scaling relation.

At each source, we remove the entire self-pole to define the regular
part, retaining its Taylor subtraction. We require its value and
gradient to vanish. Discrete scaling and polygonal symmetry reduce
these conditions to two scalar equations at one vertex; only the
radial component of the gradient survives. These are precisely the
leading equations associated with dilation and radial translation
when the poles are subsequently replaced by bubbles. Vanishing of
both terms also leaves a quadratic background near each source.

The two balance equations are solved by taking $N$ sufficiently large.
Writing $h=2\pi/N$ and $t=e^{-\alpha h}$, their solutions satisfy
\[
 \begin{cases}
 \beta=3+b_Nh,\quad b_N\longrightarrow b_*>0,&s=1,\\
 \displaystyle\beta-2=\frac{1}{\log N+O(1)},&s=2,\\
 \beta-2\asymp h^{s-2},&s>2.
 \end{cases}
\]
The limiting aspect ratio $\alpha_*$ belongs to $(2.5,2.7)$ for
$s=1$, to $(6/\pi,2)$ for $s=2$, and to $(1,\sqrt3)$ for $s>2$.
These branches reflect different balances between the nearby lattice
interaction and the circular average. They determine the geometric
parameters for the same planar configuration; the positivity argument
and the subsequent gluing construction apply to all $s\ge1$.

Positivity on the source plane follows from a common lattice limit.
For $\Lambda_\alpha=\alpha\Z\times\Z$ and $g(Y)=|Y|^{-s}$, define
\[
 D_s^\alpha(Y)=|Y|^{-s}+
 \sum_{v\in\Lambda_\alpha\setminus\{0\}}
 \bigl[g(Y-v)-g(-v)-\nabla g(-v)\cdot Y\bigr].
\]
The summands are $O(|v|^{-s-2})$ on compact sets away from poles,
so the series converges for every $s>0$. After rescaling near a
vertex, the matched point potential converges to a positive multiple
of this function. Both matching conditions are used in the
first-order subtraction: the local geometric comparison error is
$O(h\log(2/h))$ for every $s\ge1$, and the remaining terms tend to zero.
Subtracting only the value would lose the cancellation needed for
this comparison. Theta-function and Bessel-function representations
give strict positivity for the aspect ratios above. Together with
the positive pole and discrete scaling, this yields
$F(x,0)\ge c|x|^\beta$.

We next add $\kappa|z|^2$, with $\kappa$ fixed sufficiently large.
Positivity near the source plane and control on a fixed annulus give
\[
 F(X)+\kappa|z|^2\ge c\bigl(|X|^\beta+|z|^2\bigr)
 \qquad (0<|X|<1).
\]
The estimate extends to the inner annuli because $\beta>2$. This
quadratic term is biharmonic and vanishes together with its gradient
on the source plane, so it preserves the balance equations. Its
negative Laplacian, however, is the negative constant $-2(n-2)\kappa$.

Let $\lambda>0$ be small. We replace each interior pole, indexed by
$i\ge0$, with
\[
 Q_{\mu_i,q_{i\nu}}(X)=\mu_i^{-a}
 Q\bigl((X-q_{i\nu})/\mu_i\bigr),\qquad
 Q(y)=A_n(1+|y|^2)^{-a},
\]
where $\Delta^2Q=Q^p$. Matching its leading tail to the point mass gives
\[
 A_n\mu_i^a=\frac{\lambda M_i}{N},\qquad
 \mu_i\asymp\lambda^{2/s}R_i^{2+2\beta/s}\ll R_i.
\]
We subtract $T_2$ from each complete bubble shell. The exterior
point sources remain outside $B_1$ and provide a smooth biharmonic
background. Adding $\lambda\kappa|z|^2$ gives an approximate solution
$W$ satisfying $W\ge c\lambda|X|^\beta$ and, in each bubble core,
$W\ge c(Q_{\mu_i,q_{i\nu}}+\lambda|z|^2)$. The two vanishing
regular terms control the local interaction, while Taylor subtraction
controls the sum over distant shells.
Near a core, the positive bubble dominates the error caused by
smoothing. Outside the cores, the replacement error has order
$\lambda^{1+4/s}|X|^4$. Since $\beta<4$, this is smaller than the
positive planar lower bound uniformly down to the origin. Thus the
same approximation works in both regions.

We allow the shell masses and radii to vary as
$M_i e^{\eta_i}$ and $R_i e^{b_i}$, and write
$\vartheta_i=(\eta_i,b_i)$. For fixed small $\vartheta$, we solve
the projected equation for $U_\vartheta=W_\vartheta+\phi_\vartheta$.
The kernel of the single-bubble linearized operator consists of its
dilation and translation Jacobi fields. Polygonal and normal
rotational symmetries leave only dilation and radial translation
multipliers on each shell. A uniform augmented inverse, followed by
control of interactions between different cores, solves the countable
projected problem. Separate estimates for local potential amplitude
and source mass account for dimension-dependent moment integrals.
The potential amplitude controls the correction in its own core,
whereas the source mass controls its effect on other shells. Keeping
these quantities separate permits the logarithmic and power losses
in the higher-dimensional moment estimates without imposing an
incorrect uniform integrability assumption on the bubble tails.
The correction is controlled relative to $\lambda|X|^\beta$; nonlinear
estimates are made where $W_\vartheta+\phi\ge W_\vartheta/2$, which
also treats noninteger exponents and $p<2$.

To identify the remaining equations, we use the actual projected
solution and write $q_i$ for a representative center on shell $i$.
Let $f_i$ be its localized nonlinear source together with
the multiplier source in one core, let $\Gamma$ be the fundamental
solution of $\Delta^2$, and put
\[
 H_i=U_\vartheta-\Gamma*f_i,\qquad m_i=\int f_i.
\]
Then $H_i$ is biharmonic near the core and
$m_i=(\int Q^p)\mu_i^a(1+o(1))>0$.
Dilation and radial Pohozaev identities on balls of radius comparable
to $R_i$ have leading terms $a m_iH_i(q_i)$ and
$m_i\partial_rH_i(q_i)$. Self-interaction cancels by symmetry and
homogeneity of the kernel. Dividing the dilation and radial fluxes
by $a m_i\lambda R_i^\beta$ and
$m_i\lambda R_i^\beta/|q_i|$, respectively, gives
\[
 \mathcal G_\lambda(\vartheta)
   =\mathcal F_+(\vartheta)+\mathscr R_\lambda(\vartheta),\qquad
 \|\mathscr R_\lambda(0)\|_{\ell^1_\rho}
       +\Lip(\mathscr R_\lambda)\le C\lambda^{2/s}.
\]
Here $\mathcal F_+$ is the restriction to $i\ge0$ of the bilateral
point-source map with components
\[
 \mathcal F_i(\vartheta)=R_i^{-\beta}
 \bigl(F_{\vartheta,i}^{\rm reg}(q_i),
       q_i\cdot\nabla F_{\vartheta,i}^{\rm reg}(q_i)\bigr),
\]
where $F_{\vartheta,i}^{\rm reg}$ denotes the regular part of the
perturbed point field at $q_i$. We use the weighted sequence norm
$\|v\|_{\ell^1_\rho}=\sum_i\rho^i|v_i|$, with $\rho>1$ chosen below.
This comparison uses the effective backgrounds $H_i$, their source
moments and cutoff errors, rather than just the residual of $W$.
In particular, the correction contributes to the regular jet and
the localized source contributes to the mass used in normalization.
Both contributions must be included before the flux equations can
be compared with point-source balance. Estimates for the dependence
on the shell parameters give the stated Lipschitz control.

The derivative $\mathcal A=D\mathcal F(0)$ is a bilateral convolution
operator with matrix symbol $\mathbf A(\zeta)$, holomorphic in
\[
 t^{\beta-2}<|\zeta|<t^{-(4-\beta)}.
\]
We prove $\det\mathbf A(-1)>0$. Thus the determinant is not
identically zero. We can therefore choose $\rho>1$ sufficiently
close to one so that $\det\mathbf A$ does not vanish on $|\zeta|=\rho$.
The Laurent coefficients of $\mathbf A^{-1}$ give a bounded weighted
convolution inverse. For the interior equations, we retain all
interior parameters and finitely many exterior parameters. Truncating
the sufficiently distant exterior parameter columns produces a small
operator error, corrected by a Neumann series. This gives a bounded
right inverse $\mathcal R$ for the restricted derivative
$\mathcal A_+=D\mathcal F_+(0)$; the exterior background itself is retained.
The finite truncation concerns only which exterior masses and radii
are allowed to move. It does not discard the infinitely many fixed
exterior terms. This distinction preserves the balanced background
while providing a right inverse for the equations on the interior
half of the shell array.

On the range of $\mathcal R$, the reduced equations become
\[
\vartheta=-\mathcal R\bigl[
 \mathcal F_+(\vartheta)-\mathcal A_+\vartheta
            +\mathscr R_\lambda(\vartheta)\bigr].
\]
The quadratic remainder of $\mathcal F_+$ and the Lipschitz estimate
above give a contraction in a ball of radius $C\lambda^{2/s}$,
with fixed point $\vartheta_\lambda$.
The vanishing Pohozaev fluxes then eliminate both multipliers on every
shell, since their pairing with the two surviving Jacobi fields is
an invertible matrix plus a small error. Hence
$U=U_{\vartheta_\lambda}$ is an exact positive solution.
The parameter correction and the projected correction remain small
in norms that control all shells at once. In particular, positivity
does not deteriorate as the bubble centers approach the origin.

Finally, the perturbed bubble centers give
$U(q_{i\nu}(\vartheta_\lambda))\gtrsim\lambda^{-1}R_i^{-s-\beta}\to\infty$.
On the normal subspace, the solution tends to zero. Differentiated
shell estimates additionally show that the full Hessian of
$U-\lambda\kappa|z|^2$ tends to zero there, giving
\[
 \lim_{z\to0,\ z\ne0}(-\Delta U)(0,z)
       =-2(n-2)\lambda\kappa<0.
\]
This derivative conclusion uses separation of the normal subspace
from the source cores and summability of the differentiated,
compensated kernels; it does not follow merely from the pointwise
decay of the correction.
The lower estimate for the spherical minimum follows from the global
bound $U\ge c\lambda|X|^\beta$. For its upper estimate, we choose a fixed direction in
the source plane separated from all vertex directions. Along this
ray the quadratic term vanishes and the shell tails are
$O(\lambda r^\beta)$, yielding $m_U(r)\asymp r^\beta$.
The normal axis has order $r^2$; the planar ray is therefore needed
to identify this smaller spherical minimum.
The peaks, the normal subspace and this planar ray therefore give
the asserted singularity, negative Laplacian and lower envelope.

This paper is organized as follows. In Sects.~\ref{geo:section}--\ref{end:section},
we prove Theorem~\ref{intro:main}. Corollary~\ref{intro:envelope} is proved
at the end of Sect.~\ref{end:section}.
Throughout the proof, the dimension and all geometric parameters are fixed before the small
amplitude $\lambda$ is chosen. Constants may depend on these fixed
parameters, but are independent of the shell index and sufficiently
small $\lambda$. We write $f\lesssim g$ if $f\le Cg$, and
$f\asymp g$ if both inequalities hold. No estimate uniform as
$n\to\infty$ is asserted. The notation $D^j$ denotes derivatives of
order $j$; $S^{n-1}$ is the unit sphere in $\R^n$.

\section{A positive balanced configuration}\label{geo:section}

Put $s=n-4\ge1$ and write $X=(x,z)\in\R^2\times\R^{n-2}$.
For an even integer $N\ge6$, let
\[
 v_\ell=(\cos(2\pi\ell/N),\sin(2\pi\ell/N)),\quad
 P_N(X)=\frac1N\sum_{\ell=0}^{N-1}|X-(v_\ell,0)|^{-s}.
\]
The quadratic term of $P_N$ at the origin is
\begin{equation}\label{geo:H2}
 H_2(X)=c|x|^2-\frac{s}{2}|z|^2,\qquad c=\frac{s^2}{4}.
\end{equation}
Indeed, the polygon has second moment $I_2/2$. Its evenness eliminates
the odd Taylor terms. Given $0<t<1$ and $2<\beta<4$, define
\begin{equation}\label{geo:field}
 F(X)=\sum_{j\in\Z}t^{j\beta}
 \bigl[P_N(t^{-j}X)-1-t^{-2j}H_2(X)\bigr].
\end{equation}
This series converges locally smoothly away from the origin and the points
$q_{j,\ell}=t^j(v_\ell,0)$. The outer shells have size
$O(t^{j(\beta-4)})$ as $j\to-\infty$, and the inner quadratic terms have
size $O(t^{j(\beta-2)})$ as $j\to+\infty$. At $q_{j,\ell}$, the
\emph{regular part} means that the entire term
$N^{-1}t^{j(s+\beta)}|X-q_{j,\ell}|^{-s}$ has been removed.

\begin{proposition}\label{geo:seed}
For each integer $s\ge1$, all sufficiently large even integers $N$ admit
parameters $t\in(0,1)$ and $\beta\in(2,4)$, with $3<\beta<4$ when
$s=1$ and $2<\beta<3$ when $s\ge2$, such that the regular value
and gradient of $F$ vanish at every $q_{j,\ell}$. Moreover,
\begin{equation}\label{geo:floor}
 F(x,0)\ge c_0|x|^\beta
\end{equation}
away from the sources, where $c_0>0$. There are $\kappa,c_1>0$ such that
\begin{equation}\label{geo:ambient-positive}
 F(X)+\kappa|z|^2\ge c_1(|X|^\beta+|z|^2),\qquad 0<|X|\le1.
\end{equation}
The constants may depend on $s$ and on the chosen $N$.
\end{proposition}

We prove the proposition below. All limits in this section keep $s$ fixed.
Set $h=2\pi/N$, $T=\alpha h$, and $t=e^{-T}$. Denote the regular value
and radial derivative at $(v_0,0)$ by $V_N$ and $R_N$, respectively.
The dihedral symmetry and the $O(n-2)$ symmetry eliminate all other
components of the gradient. Thus the required equations are $V_N=R_N=0$.

\subsection{Positive coefficients of the polygon kernel}

\begin{lemma}\label{geo:coeff}
Let $b_m=(s/2)_m/m!$, where $(s/2)_m$ is the rising factorial. On the
ray through $v_0$, for $0\le r<1$,
\begin{align}
 P_N(rv_0,0)&=\sum_{\substack{d\ge0\\d\ \mathrm{even}}}C_{d,N}r^d,
 \label{geo:coeff-formula}\\
 C_{d,N}&=\sum_{\substack{k\in\Z\\|k|\le d/N}}
 b_{(d+kN)/2}b_{(d-kN)/2}>0.\notag
\end{align}
In particular, $C_{0,N}=1$ and $C_{2,N}=c$. Put
\begin{equation}
 S_N=\frac1N\sum_{\ell=1}^{N-1}(2\sin(\pi\ell/N))^{-s}.\label{geo:SN}
\end{equation}
For $s\ge2$, define
\begin{align}
 W_s(y)&=\frac1{2^{s-2}\Gamma(s/2)^2}
       \sum_{|k|<y}(y^2-k^2)^{s/2-1},\qquad y>0,\label{geo:Ws}\\
 G_s(u)&=\sum_{k\in\Z}(u^2+k^2)^{-s/2},\qquad u>0.\label{geo:Gs}
\end{align}
For $s\ge2$ and noninteger $y>0$, as $d/N\to y$ through even $d$,
\begin{equation}\label{geo:coeff-limit}
 N^{2-s}C_{d,N}\longrightarrow W_s(y),\qquad
 N^{2-s}C_{d,N}\le C_s(y+N^{-1})^{s-2}(1+y).
\end{equation}
For $s\ge2$, we also have
\begin{equation}\label{geo:laplace}
 h^{s-1}S_N\longrightarrow\frac{2\zeta(s)}{2\pi},\qquad
 G_s(u)=\frac{(2\pi)^s}{2}\int_0^\infty W_s(y)e^{-2\pi uy}\dd y.
\end{equation}
\end{lemma}

\begin{proof}
Expand $(1-re^{i\theta})^{-s/2}(1-re^{-i\theta})^{-s/2}$ and average
over $\theta=2\pi\ell/N$. The average retains precisely the indices
whose difference belongs to $N\Z$, giving \eqref{geo:coeff-formula}.
Stirling's formula and the estimate
$b_m\le C_s(m+1)^{s/2-1}$ give \eqref{geo:coeff-limit}.
The first limit in \eqref{geo:laplace} follows from
$\sin(\pi\ell/N)\ge2\ell/N$ for $1\le\ell\le N/2$ and dominated
convergence. For the second identity, apply
\[
 r^{-s}=\frac1{\Gamma(s/2)}\int_0^\infty q^{s/2-1}e^{-qr^2}\dd q
\]
and Poisson summation to the $k$-sum defining $G_s$. Equivalently, the
integral of each $(y^2-k^2)_+^{s/2-1}$ in \eqref{geo:Ws} is the standard
positive Laplace integral for the modified Bessel function
$K_{(s-1)/2}$. The duplication formula for $\Gamma$ gives the factor
$(2\pi)^s/2$. All integrands are nonnegative, so the exchanges of sums
and integrals are justified by monotone convergence.
\end{proof}

In every subsequent sum over $d$, the index is even. Write
$H(u)=(e^u-1)^{-1}$. Summing first over the shells gives the exact formulas
\begin{align}
 V_N={}&S_N-1-c+\sum_{d\ge0}C_{d,N}H(T(s+\beta+d))
       -H(T\beta)-cH(T(\beta-2))\notag\\
     &+\sum_{d\ge4}C_{d,N}H(T(d-\beta)),\label{geo:Vexact}\\
 R_N={}&-\frac{s}{2}S_N-2c
       -\sum_{d\ge0}(s+d)C_{d,N}H(T(s+\beta+d))\notag\\
     &-2cH(T(\beta-2))
       +\sum_{d\ge4}dC_{d,N}H(T(d-\beta)).\label{geo:Rexact}
\end{align}
For fixed $N,T,\beta$, each displayed series is absolutely convergent.
These identities also follow directly from homogeneity by pairing the
shells $t^j$ and $t^{-j}$.

\subsection{The lattice pressure}

For $s>2$, define
\begin{equation}\label{geo:Epstein}
 E_s(\alpha)=\sum_{(j,k)\ne(0,0)}(\alpha^2j^2+k^2)^{-s/2},
 \qquad g_s(\alpha)=(\alpha E_s(\alpha))'.
\end{equation}

\begin{lemma}\label{geo:pressure}
For every integer $s\ge3$, $g_s$ has exactly one zero
$\alpha_s\in(1,\sqrt3)$, and $g_s'(\alpha_s)>0$.
In addition, $(u^{s-1}G_s(u))'<0$ for $u>0$.
\end{lemma}

\begin{proof}
Interchanging the two lattice coordinates gives
$g_s(1)=(1-s/2)E_s(1)<0$. Put $\nu=(s-1)/2$. Poisson summation gives
\begin{align}
 E_s(\alpha)={}&2\zeta(s)+
 \frac{2\sqrt\pi\Gamma(\nu)}{\Gamma(s/2)}
       \zeta(s-1)\alpha^{1-s}\notag\\
 &+\frac{8\pi^{s/2}}{\Gamma(s/2)}\alpha^{-\nu}
       \sum_{j,k\ge1}(k/j)^\nu K_\nu(2\pi\alpha jk).
 \label{geo:poisson}
\end{align}
The integral $K_\nu(u)=\int_0^\infty e^{-u\cosh v}\cosh(\nu v)\dd v$
shows that $K_\nu$ is completely monotone. In $\alpha E_s(\alpha)$,
the nonconstant terms other than the linear term are positive multiples
of $\alpha^{2-s}$ and $\alpha^{(3-s)/2}K_\nu(c\alpha)$.
Their second derivatives are nonnegative, and that of
$\alpha^{2-s}$ is strictly positive. Thus $g_s'>0$.

It remains to check $g_s(\sqrt3)>0$. For $s\ge5$, the absolute value
of the negative contribution is bounded by
\begin{align*}
 N_s={}&2(s-1)3^{-s/2}\zeta(s)\\
 &+2(s-1)\frac{\sqrt\pi\Gamma((s+1)/2)}{\Gamma((s+2)/2)}
                  3^{(1-s)/2}\zeta(s-1).
\end{align*}
Indeed, discard the positive $k^2$ in the numerator and bound the
sum over $k\in\Z$ by its central term plus the integral over $\R$.
The resulting summand
$3(s-1)j^2(3j^2+k^2)^{-s/2-1}$ decreases with $s\ge3$.
Consequently the negative contribution for $s\ge5$ is at most
\[
 N_5=\frac8{9\sqrt3}\zeta(5)+\frac{128}{135}\zeta(4)<1.6,
\]
whereas the $j=0$ contribution is $2\zeta(s)>2$.
For $s=3,4$, the first two terms of \eqref{geo:poisson} give, respectively,
\[
 2\zeta(3)-\frac43\zeta(2)>0.2,\qquad
 2\zeta(4)-\frac{2\pi}{3\sqrt3}\zeta(3)>0.69.
\]
At $\alpha=\sqrt3$, put $q=e^{-c_*}$ and $c_*=2\pi\alpha$.
The absolute values of the derivatives of the remaining terms in
$\alpha E_s$ are bounded by
\begin{align*}
 T_3&\le8\pi\alpha^{-3/2}\zeta(2)(c_*+1+c_*^{-1})
                 \frac{q(1+q)}{(1-q)^3}<0.005,\\
 T_4&\le4\pi^2\alpha^{-2}\zeta(3)(c_*+2+2c_*^{-1})
                 \frac{q(1+q)}{(1-q)^3}<0.005.
\end{align*}
These bounds follow from $K_0\le K_{1/2}$, $K_1\le K_{3/2}$,
the explicit half-integer formulas, and
$\sum_{m\ge1}m^2q^m=q(1+q)/(1-q)^3$.
The displayed numerical inequalities follow, for example, by bounding
each zeta tail by its defining integral. They leave a strict margin.
This proves the assertion about $g_s$.
Finally, the one-row Poisson formula for $G_s$, together with
$(u^\nu K_\nu(cu))'=-cu^\nu K_{\nu-1}(cu)<0$, proves the last assertion.
\end{proof}

\subsection{Positivity of the lattice potential}

For $\Lambda_\alpha=\alpha\Z\times\Z$ and $g(w)=|w|^{-s}$, define
\begin{equation}\label{geo:lattice-difference}
 D_s^\alpha(Y)=|Y|^{-s}+
 \sum_{v\in\Lambda_\alpha\setminus\{0\}}
 [g(Y-v)-g(-v)-\nabla g(-v)\cdot Y].
\end{equation}
The summands are $O_K(|v|^{-s-2})$ on compact sets off the lattice,
so the series is absolutely and locally uniformly convergent for
every $s>0$. Opposite vectors cancel the linear terms. For $s>2$,
it equals the lattice potential minus $E_s(\alpha)$.
Put $\Theta_\alpha(u;Y)=\sum_{v\in\Lambda_\alpha}e^{-u|Y-v|^2}$.
Then
\begin{equation}\label{geo:lattice-mellin}
 D_s^\alpha(Y)=\frac1{\Gamma(s/2)}\int_0^\infty
 u^{s/2-1}[\Theta_\alpha(u;Y)-\Theta_\alpha(u;0)+1]\dd u.
\end{equation}
To justify the interchange, pair $v$ with $-v$. For large $|v|$,
the Gaussian second difference is bounded by
$C_K(u+u^2|v|^2)e^{-cu|v|^2}$; its Mellin integral is
$O_K(|v|^{-s-2})$. This also proves periodicity of $D_s^\alpha$.

\begin{lemma}\label{geo:lattice-positive}
For $s>0$ and $1\le\alpha\le\sqrt3$, the minimum of
$D_s^\alpha$ is attained at $(\alpha/2,1/2)$ and is strictly positive.
For fixed $s>0$, these minima have a positive lower bound on
$[1,\sqrt3]$. For $s=2$, the same minimum property holds, and
\begin{equation}\label{geo:critical-gap}
 \min D_2^\alpha>0.99,\qquad \sqrt3\le\alpha\le2.
\end{equation}
For $s=1$, one has
\begin{equation}\label{geo:subcritical-gap}
 \min D_1^\alpha>0.53,\qquad 5/2\le\alpha\le27/10.
\end{equation}
\end{lemma}

\begin{proof}
The periodic Gaussian sum on a line is minimized halfway between
successive lattice points. To see this, apply Poisson summation and
the Jacobi product formula: each factor depending on the shift has the
form $1+2q^{2m-1}\cos(2\pi y)+q^{4m-2}$ and is minimized at
$y=1/2$. Multiply the two positive one-dimensional sums and integrate
against the positive Mellin weight in \eqref{geo:lattice-mellin}.
This proves the minimum property for every $s>0$.

Write
\[
 \vartheta_3(u)=\sum_{k\in\Z}e^{-uk^2},\quad
 \vartheta_2(u)=\sum_{k\in\Z}e^{-u(k+1/2)^2},\quad
 \vartheta_4(u)=\sum_{k\in\Z}(-1)^ke^{-uk^2}.
\]
The duplication and cubic identities are
\begin{align}
 \vartheta_3(u)^2-\vartheta_2(u)^2&=\vartheta_4(u/2)^2,\notag\\
 \vartheta_3(u)\vartheta_3(3u)-\vartheta_2(u)\vartheta_2(3u)
 &=\vartheta_4(u)\vartheta_4(3u).\label{geo:theta-identities}
\end{align}
Both follow by changes of lattice indices. For the second, the map
$(m,k)\mapsto((m+3k)/2,(m-k)/2)$ preserves $m^2+3k^2$ and maps
the class $m+k$ odd to one of the two equal half-integer classes;
reflection exchanges those classes.
The product formula shows that $\vartheta_4(u/2)^2$ is strictly increasing.
Writing $f=\vartheta_3$, $g=\vartheta_2$, the first identity gives
$(-f')/(-g')<g/f$, and $g/f$ decreases.
Hence, for $v\ge u$,
\[
 \partial_v[f(u)f(v)-g(u)g(v)]
 =(-g'(v))\left[g(u)-\frac{-f'(v)}{-g'(v)}f(u)\right]>0.
\]
If $u\le v\le3u$, \eqref{geo:theta-identities} therefore gives
\[
 f(u)f(v)-g(u)g(v)\le\vartheta_4(u)\vartheta_4(3u)<1.
\]
The value of $D_s^\alpha$ at the cell center is
\begin{equation}\label{geo:mellin-gap}
 \frac1{\Gamma(s/2)}\int_0^\infty u^{s/2-1}
 [1-f(u)f(\alpha^2u)+g(u)g(\alpha^2u)]\dd u.
\end{equation}
For $1\le\alpha\le\sqrt3$, the integrand is strictly positive.
Continuity and compactness give the uniform positive lower bound.

For $s=2$, summation in the second coordinate gives
\begin{align}
 D_2^\alpha(\alpha/2,1/2)
 ={}&\frac{4\pi\log2}{\alpha}-\frac{\pi^2}{3}\notag\\
 &+\frac{2\pi}{\alpha}\sum_{j\ge0}
 \frac{\tanh(\pi\alpha(j+1/2))-1}{j+1/2}\notag\\
 &-\frac{2\pi}{\alpha}\sum_{j\ge1}
 \frac{\coth(\pi\alpha j)-1}{j}.\label{geo:critical-gap-series}
\end{align}
Here one uses
$\sum_k(k^2+b^2)^{-1}=\pi\coth(\pi b)/b$ and the corresponding
half-integer formula $\pi\tanh(\pi b)/b$.
The two harmonic divergences are subtracted with the same cutoffs;
their difference is $4\pi\log2/\alpha$.
For $\sqrt3\le\alpha\le2$, the first line is greater than $1.06$.
With $q=e^{-\pi\alpha}<0.0044$, the absolute values of the next two
lines are at most
\[
 \frac{8\pi}{\alpha}\operatorname{arctanh}q<0.064,
 \qquad
 \frac{4\pi}{\alpha}\frac{-\log(1-q^2)}{1-q^2}<0.00015.
\]
This proves \eqref{geo:critical-gap}.

For $s=1$, let $\gamma$ be Euler's constant. One-dimensional Poisson summation yields
\begin{align}
 D_1^\alpha(\alpha/2,1/2)
 ={}&2\log(2\pi/\alpha)-2\gamma\notag\\
 &+8\sum_{m\ge1}\left[
 (-1)^m\sum_{j\ge0}K_0(2\pi\alpha m(j+1/2))
       -\sum_{j\ge1}K_0(2\pi\alpha mj)\right].
 \label{geo:s1-gap-series}
\end{align}
For clarity, derive the identity first for $s>2$. The zero Fourier
mode of the center difference is
\[
 -2\zeta(s)+2C_s\alpha^{1-s}(2^{s-1}-2)\zeta(s-1),
 \qquad C_s=\frac{\sqrt\pi\Gamma((s-1)/2)}{\Gamma(s/2)}.
\]
The nonzero modes converge exponentially. The difference definition
\eqref{geo:lattice-difference} is analytic for $\operatorname{Re}s>0$,
and the right side continues there after its removable singularities
are canceled. At $s=1+\varepsilon$,
\[
 C_s=2/\varepsilon+2\log2+O(\varepsilon),\qquad
 (2^\varepsilon-2)\zeta(\varepsilon)
       =\tfrac12+\tfrac12\log\pi\,\varepsilon+O(\varepsilon^2).
\]
The poles cancel and give the first line of \eqref{geo:s1-gap-series};
the nonzero-mode factor tends to $4K_0$, giving the second line.
With $q=e^{-\pi\alpha}$ and
$K_0(u)\le\sqrt{\pi/(2u)}e^{-u}$, the absolute value of that line is at most
\[
 \mathcal E(\alpha)=\frac8{\sqrt{2\alpha}}
 \frac{q}{(1-q)(1-q^2)}
 +\frac4{\sqrt\alpha}\frac{q^2}{(1-q^2)^2}.
\]
Each factor decreases in $\alpha$. For $5/2\le\alpha\le27/10$,
the first line exceeds $0.5348$, whereas $\mathcal E(5/2)<0.0014$.
This proves \eqref{geo:subcritical-gap}.
\end{proof}

\subsection[The matching equations for s > 2]{The matching equations for $s>2$}

\begin{lemma}\label{geo:noncritical-root}
For every integer $s\ge3$ and all sufficiently large even $N$, the
equations $V_N=R_N=0$ have a solution with
\[
 \alpha_N\longrightarrow\alpha_s,\qquad
 \beta_N=2+b_Nh^{s-2},\qquad
 b_N\longrightarrow\frac{2\pi c}{\alpha_sE_s(\alpha_s)}>0.
\]
\end{lemma}

\begin{proof}
Put $\beta=2+bh^{s-2}$, where $(\alpha,b)$ ranges over a compact
subset of $(0,\infty)^2$. Lemma~\ref{geo:coeff} gives the following
limits in $C^1$ on such compact sets:
\begin{align}
 h^{s-1}V_N&\longrightarrow
 \frac{E_s(\alpha)}{2\pi}-\frac{c}{\alpha b},\label{geo:Vlimit}\\
 h^{s-1}R_N&\longrightarrow
 \frac{-sE_s(\alpha)-(s+4)\alpha E_s'(\alpha)}{4\pi}
       -\frac{2c}{\alpha b}.\label{geo:Rlimit}
\end{align}
We provide the domination for these assertions. The only slow term is
$H(T(\beta-2))$, which is kept separate and has the displayed explicit
$C^1$ limit. For $d\ge4$, set $y=d/N$. The coefficient bound in
\eqref{geo:coeff-limit}, multiplied by the Bose kernel $H$, has near-zero
majorant $Cy^{s-3}$ and an exponential majorant at infinity.
To combine the two differentiated sums in \eqref{geo:Rexact}, use
\begin{align}
 &d\{H(T(d-\beta))-H(T(d+s+\beta))\}\notag\\
 &\qquad=(s+2\beta)dT\int_0^1
       -H'(T(d-\beta+\theta(s+2\beta)))\dd\theta.\label{geo:difference-quotient}
\end{align}
The combined limiting kernel is
$-sH(2\pi\alpha y)-(s+4)(2\pi\alpha y)H'(2\pi\alpha y)$.
For its $\alpha$ derivative the same bound follows from
$yH'=O(y^{-1})$ and $y^2H''=O(y^{-1})$ near zero.
For the $\beta$ derivative, before taking the limit, the bound is
$C(y^{-1}+N^{-1}y^{-2})\le C'y^{-1}$, since $Ny=d\ge4$.
The $b$ derivative has the additional factor $h^{s-2}$.
The finitely many $d=0,2$ terms and their relevant derivatives are
$O(h^{-1})$ and disappear after normalization.
Splitting the Riemann sums into $(0,\varepsilon)$, a compact interval,
and the exponential tail now proves $C^1$ convergence.
The Laplace identity \eqref{geo:laplace} identifies the limits as
\eqref{geo:Vlimit}--\eqref{geo:Rlimit}.

The first limiting equation determines $b=2\pi c/(\alpha E_s(\alpha))$.
On this curve, the second is
$-(s+4)(E_s+\alpha E_s')/(4\pi)=0$.
Lemma~\ref{geo:pressure} gives a simple zero, and the derivative of the
first equation in $b$ is positive. The limiting two-dimensional zero
is therefore nondegenerate. Stability of a nondegenerate zero under
$C^1$ perturbations proves the lemma.
\end{proof}

\subsection[The logarithmic case s = 2]{The logarithmic case $s=2$}

\begin{lemma}\label{geo:critical-root}
For $s=2$ and all sufficiently large even $N$, the equations
$V_N=R_N=0$ have a solution with
\[
 \alpha_N\longrightarrow\alpha_2\in(6/\pi,2),\qquad
 \beta_N=2+\frac1{\log N+O(1)}.
\]
Here $\alpha_2$ is the unique zero of
\begin{equation}\label{geo:D2}
 D_2(\alpha)=-\frac1{12}+\frac1{2\pi\alpha}
       +2\sum_{j\ge1}\frac{e^{2\pi\alpha j}}{(e^{2\pi\alpha j}-1)^2}.
\end{equation}
\end{lemma}

\begin{proof}
In polar coordinates on the source plane, geometric summation gives
\begin{equation}\label{geo:critical-P}
 P_N(r,\theta)=\frac{1-r^{2N}}
 {(1-r^2)(1-2r^N\cos(N\theta)+r^{2N})},\qquad 0<r<1,
\end{equation}
and $S_N=(N^2-1)/(12N)$.
Write $\delta=\beta-2$ and $b=3+\delta$. Pairing the shells yields
\begin{align}
 V_N={}&S_N-2+\sum_{j\ge1}\left[
 \frac{\cosh(bjT)}{\sinh(jT)}\coth(NjT/2)
       -2\cosh(\beta jT)-2\cosh(\delta jT)\right],\label{geo:critical-V}\\
 R_N-\beta V_N={}&2b-4-bS_N-\sum_{j\ge1}G_{N,\delta}'(jT),
 \label{geo:critical-R}\\
 G_{N,\delta}(u)={}&\frac{\sinh(bu)}{\sinh u}\coth(Nu/2)
                 -2\sinh(\beta u)-2\sinh(\delta u).\notag
\end{align}
All cancellations in these expressions are taken within each bracket.
Uniformly for $\alpha$ in a fixed compact positive interval and for
$\delta$ near $1/\log N$,
\begin{align}
 \frac{V_N}{N}&=\frac{\log N-\delta^{-1}}{2\pi\alpha}+O(1),
 \label{geo:critical-Vasymp}\\
 \frac1N\partial_\delta V_N
 &=\frac1{2\pi\alpha\delta^2}+O(1).\label{geo:critical-Vderiv}
\end{align}
Indeed, remove the factor $\coth(Nu/2)-1$ in
\eqref{geo:critical-V}; its contribution, divided by $N$, and its
$\delta$ derivative are bounded by a convergent exponential sum in $j$.
In the remaining bracket, addition of $e^{-\delta u}$ leaves a function
with expansion $u^{-1}+O(1)$ near zero and uniform exponential decay at
infinity. Its $\delta$ derivative is bounded near zero and integrable
at infinity. Thus its Riemann sum is $T^{-1}\log(T^{-1})+O(T^{-1})$,
whereas $\sum_{j\ge1}e^{-\delta jT}=H(T\delta)$ supplies the slow term.
This proves \eqref{geo:critical-Vasymp}--\eqref{geo:critical-Vderiv}.
Consequently $V_N=0$ has a continuous branch
$\delta_N(\alpha)=(\log N+O(1))^{-1}$.

For the force equation write $G_{N,\delta}=G_{0,\delta}+H_{N,\delta}$,
where $H_{N,\delta}=(\sinh(bu)/\sinh u)(\coth(Nu/2)-1)$.
The exact identity
\begin{equation}\label{geo:G0}
 G_{0,\delta}(u)=e^{-\delta u}+e^{-(2-\delta)u}
       +e^{-(2+\delta)u}
       +e^{-3u}\frac{\sinh(\delta u)}{\sinh u}
\end{equation}
gives $\|G_{0,\delta}''\|_{L^1(0,\infty)}\le C$ and
$\int_0^\infty-G_{0,\delta}'=3+\delta=b$.
The elementary sum--integral bound for a function with integrable
derivative therefore gives
\[
 \sum_{j\ge1}-G_{0,\delta}'(jT)=\frac bT+O(1).
\]
This retains the contribution of $e^{-\delta u}$ even as $\delta\to0$.
For the other term, direct differentiation at $u=jT$ gives
\[
 \frac1{bN}\sum_{j\ge1}-H_{N,\delta}'(jT)
 \longrightarrow\frac12\sum_{j\ge1}\operatorname{csch}^2(\pi\alpha j).
\]
The summands and their errors are bounded by $C(1+j)^2e^{-c j}$,
uniformly in the stated parameter range. Equation~\eqref{geo:critical-R}
therefore gives, on $V_N=0$,
\[
 \frac{R_N}{bN}\longrightarrow D_2(\alpha)
\]
uniformly. Every summand in \eqref{geo:D2} is strictly decreasing in
$\alpha$. Also $D_2(6/\pi)>0$, while
$D_2(2)<0$: the sum there is less than $0.000008$ and
$1/12-1/(4\pi)>0.0037$. The intermediate value theorem gives the
required exact roots, and their accumulation points must be the
unique zero $\alpha_2$.
\end{proof}

\subsection[The integrable case s = 1]{The integrable case $s=1$}

\begin{lemma}\label{geo:subcritical-root}
For $s=1$ and every sufficiently large even $N$, the equations
$V_N=R_N=0$ have a solution satisfying
\[
 \alpha_N\longrightarrow\alpha_1\in(5/2,27/10),\qquad
 \beta_N=3+b_Nh,\qquad
 b_N\longrightarrow-\frac{9\alpha_1E_1(\alpha_1)}{8\pi}>0,
\]
where
\begin{equation}\label{geo:E1}
 E_1(\alpha)=2\log\alpha+2\gamma-2\log(4\pi)
             +8\sum_{j,k\ge1}K_0(2\pi\alpha jk)
\end{equation}
and $\alpha_1$ is the unique zero of $g_1=(\alpha E_1)'$ in
$[5/2,27/10]$. Here $\gamma$ is Euler's constant.
\end{lemma}

\begin{proof}
Put $c=1/4$, $\sigma=\beta+1/2$, and define
\begin{align*}
 K_N(u)&=\frac1N\sum_{\ell=0}^{N-1}
   [2(\cosh u-\cos(2\pi\ell/N))]^{-1/2},\\
 K(u)&=\frac1{2\pi}\int_0^{2\pi}
   [2(\cosh u-\cos\theta)]^{-1/2}\dd\theta,
 &\delta_N(u)&=K_N(u)-K(u),\\
 f_\beta(u)&=2\cosh(\sigma u)K(u)-2\cosh(\beta u)
                    -2c\cosh((\beta-2)u),\\
 G_\beta(u)&=2\sinh(\sigma u)K(u)-2\sinh(\beta u)
                    -2c\sinh((\beta-2)u).
\end{align*}
All cancellations are taken within these complete expressions.
Since $P_N(e^uv_0,0)=e^{-u/2}K_N(u)$, pairing opposite shells gives
\begin{align}
 V_N={}&S_N-1-c+\sum_{j\ge1}f_\beta(jT)
           +2\sum_{j\ge1}\cosh(\sigma jT)\delta_N(jT),
             \label{geo:s1-pair-value}\\
 R_N-\beta V_N={}&\beta+c(\beta-2)-\sigma S_N
       -\sum_{j\ge1}G_\beta'(jT)
       -2\sum_{j\ge1}[\sinh(\sigma u)\delta_N(u)]'_{u=jT}.
             \label{geo:s1-pair-force}
\end{align}

We first consider the angular average. Define the complete elliptic integral
$\mathcal K(k)=\int_0^{\pi/2}(1-k^2\sin^2v)^{-1/2}\dd v$. Then
\[
 K(u)=\frac{\mathcal K(\operatorname{sech}(u/2))}{\pi\cosh(u/2)}
     =\frac1\pi\log\frac8u+O(u^2(1+|\log u|)).
\]
The error may be differentiated twice, with respective bounds
$O(u(1+|\log u|))$ and $O(1+|\log u|)$. Thus, as $u\downarrow0$,
\begin{align*}
 f_\beta(u)&=-\frac2\pi\log u+\frac2\pi\log8-2(1+c)
                       +O(u^2(1+|\log u|)),\\
 G_\beta'(u)&=-\frac{2\sigma}\pi\log u
   +\frac{2\sigma}\pi(\log8-1)-2\beta-2c(\beta-2)
                       +O(u^2(1+|\log u|)).
\end{align*}
These estimates and their parameter derivatives are uniform for
$\beta$ in compact subsets of $(2,4)$. At infinity the complete
expressions and their derivatives decay exponentially. Moreover,
$G_\beta(0)=G_\beta(\infty)=0$, so $\int_0^\infty G_\beta'=0$.

For an integrable function $f(u)=-A\log u+B+O(u(1+|\log u|))$, with an
integrable derivative after subtraction of the cutoff logarithmic
term, logarithmic Euler summation gives
\begin{equation}\label{geo:log-euler}
 \sum_{j\ge1}f(jT)=\frac1T\int_0^\infty f(u)\dd u
          +\frac A2\log\frac{T}{2\pi}-\frac B2+o(1).
\end{equation}
Indeed, sum $\log(jT)$ explicitly on a truncated interval and apply
Stirling's formula. Ordinary trapezoidal summation on its complement
cancels the artificial endpoint terms. For the remainder, the
oscillatory Bernoulli integral tends to zero because its derivative
belongs to $L^1$. Here the local bounds and uniform exponential tails
make the remainder derivatives continuous in $L^1$ as the parameters
vary. Their image over a compact parameter set is compact in $L^1$,
so the error is uniform. The same argument applies to parameter
derivatives and, when needed, to $uf'(u)$.

Subtracting the two endpoint singularities from the cosecant sum gives
\[
 S_N=\frac1\pi\left(\log\frac4h+\gamma\right)+o(1).
\]
More explicitly, subtract $(2\pi x)^{-1}+[2\pi(1-x)]^{-1}$ from
$(2\sin\pi x)^{-1}$. The remainder is continuous on $[0,1]$, with
integral $\pi^{-1}\log(2/\pi)$; its Riemann sum together with
$\pi^{-1}\sum_{\ell=1}^{N-1}\ell^{-1}$ gives the formula.
Write $\mathcal A(\beta)=\int_0^\infty f_\beta(u)\dd u$.
Applying \eqref{geo:log-euler} to
\eqref{geo:s1-pair-value}--\eqref{geo:s1-pair-force} shows that
the parts without $\delta_N$ are, respectively,
\[
 \frac{\mathcal A(\beta)}{\alpha h}
       +\frac{\gamma+\log\alpha-\log(4\pi)}\pi+o(1),
 \qquad
 -\frac\sigma\pi[\gamma+\log\alpha-\log(4\pi)+1]+o(1).
\]

We next estimate the nonzero angular modes. For
$b_r=(1/2)_r/r!$, the $m$th complex Fourier coefficient of the
angular kernel $[2(\cosh u-\cos\theta)]^{-1/2}$ is
\[
 a_m(u)=e^{-(m+1/2)u}\sum_{r\ge0}b_{r+m}b_re^{-2ru},
 \qquad \delta_N(u)=2\sum_{k\ge1}a_{kN}(u).
\]
Since $b_r\le C(r+1)^{-1/2}$ and $b_r\sim(\pi r)^{-1/2}$,
Riemann summation with $r=mx$ gives, for fixed $\lambda>0$,
\begin{equation}\label{geo:s1-fourier-limit}
 a_m(\lambda/m)\longrightarrow
 \frac{e^{-\lambda}}\pi\int_0^\infty
       \frac{e^{-2\lambda x}}{\sqrt{x(x+1)}}\dd x
       =\frac1\pi K_0(\lambda).
\end{equation}
The substitution $2x+1=\cosh v$ gives the last equality. The
integrable endpoint singularity and exponential tail also justify
two derivatives in $\lambda$. For $L=0,1,2$ and $mu\ge\lambda_0>0$,
\begin{equation}\label{geo:s1-fourier-bound}
 |(u\partial_u)^L a_m(u)|
       \le C_{L,\lambda_0}(1+mu)^{L+1}e^{-mu/2}.
\end{equation}
For $u\le1$, this follows from $b_{r+m}\le Cm^{-1/2}$ and
$\sum_{r\ge0}(r+1)^{-1/2}(1+ru)^Le^{-2ru}\le C_Lu^{-1/2}$;
for $u\ge1$, separate $r=0$ from the exponential tail.
At $m=kN$ and $u=jT$, one has $mu=2\pi\alpha jk$.
For large $N$, the hyperbolic factors in
\eqref{geo:s1-pair-value}--\eqref{geo:s1-pair-force} preserve a
summable bound $C(1+jk)^C e^{-cjk}$, uniformly on compact parameter
sets under consideration. Thus \eqref{geo:s1-fourier-limit} gives
\begin{align*}
 2\sum_{j\ge1}\cosh(\sigma jT)\delta_N(jT)
     &\longrightarrow\frac4\pi\sum_{j,k\ge1}K_0(z),\\
 -2\sum_{j\ge1}[\sinh(\sigma u)\delta_N(u)]'_{u=jT}
     &\longrightarrow-\frac{4\sigma}\pi
                \sum_{j,k\ge1}[K_0(z)-zK_1(z)],
 \qquad z=2\pi\alpha jk.
\end{align*}
Parameter derivatives have the same domination. Combining the two
parts proves the uniform asymptotics
\begin{align}
 V_N&=\frac{\mathcal A(\beta)}{\alpha h}
                      +\frac{E_1(\alpha)}{2\pi}+o(1),
                    \label{geo:s1-value-limit}\\
 R_N-\beta V_N&=-\frac{\beta+1/2}{2\pi}g_1(\alpha)+o(1).
                    \label{geo:s1-force-limit}
\end{align}

The continuous coefficient is
\begin{equation}\label{geo:s1-continuous}
 \mathcal A(\beta)=
 \frac{\Gamma(-\beta/2)\Gamma((\beta+1)/2)}
 {2\Gamma((\beta+2)/2)\Gamma((1-\beta)/2)},
 \qquad \mathcal A(3)=0,\quad\mathcal A'(3)=\frac49.
\end{equation}
For a direct derivation, first take $-1<\beta<0$ and compute
$(2\pi)^{-1}\int_{\R^2}|e_1-y|^{-1}|y|^{\beta-1}\dd y$.
Represent the two powers by Gaussian integrals, integrate in $y$,
and set the two Gaussian parameters equal to $\lambda a$ and
$\lambda(1-a)$. The $\lambda$ integral gives
$\Gamma(-\beta/2)$, and the $a$ integral gives
$B(1/2,(\beta+1)/2)$, yielding the quotient in
\eqref{geo:s1-continuous}. Split the radial integral at one, subtract
$1+cq^{-2}$ on $q>1$, and add $-1/\beta-c/(\beta-2)$.
This continues the same expression to $-1<\operatorname{Re}\beta<4$
away from its poles. For $2<\beta<4$, it equals the integral with
the two subtractions on all $q>0$, namely $\mathcal A(\beta)$.
The residue of $\Gamma$ at $-1$ gives its value and derivative at three.

Set $\beta=3+bh$. Uniformly for
$(\alpha,b)\in[5/2,27/10]\times[1,3]$, we obtain
\begin{equation}\label{geo:s1-rectangle-limit}
 V_N\longrightarrow\frac{4b}{9\alpha}+\frac{E_1(\alpha)}{2\pi},
 \qquad R_N-\beta V_N\longrightarrow-\frac7{4\pi}g_1(\alpha).
\end{equation}
We record strict bounds that determine the signs. On this interval,
\[
 -2.076<E_1(\alpha)<-1.920,\qquad
 g_1(5/2)<-0.075,\quad g_1(27/10)>0.078,\quad g_1'>0.73.
\]
These follow from $K_0\le K_{1/2}$, $K_1\le K_{3/2}$ and the
explicit half-integer formulas. Indeed, put $c_*=5\pi$,
$q=e^{-c_*}$, $A_*=(\pi/(2c_*))^{1/2}$, and
$S_r=\sum_{m\ge1}m^rq^m$. The Bessel contributions to $E_1$, $g_1$
and $g_1'$ have absolute values at most
\[
 8A_*S_1<3.82\cdot10^{-7},\quad
 8A_*(c_*S_2+2S_1)<6.76\cdot10^{-6},\quad
 \frac{16A_*}{5}(c_*^2S_3+c_*S_2+S_1)<4.02\cdot10^{-5}.
\]
Here the number of pairs with $jk=m$ is at most $m$;
$S_1=q/(1-q)^2$, $S_2=q(1+q)/(1-q)^3$, and
$S_3=q(1+4q+q^2)/(1-q)^4$ give explicit bounds.
The elementary parts of $g_1$ at the two endpoints are less than
$-0.07503$ and greater than $0.07888$, respectively.

The first limit in \eqref{geo:s1-rectangle-limit} is negative at
$b=1$ and positive at $b=3$. The second is positive at $\alpha=5/2$
and negative at $\alpha=27/10$. Uniform convergence preserves these
boundary signs. The rectangular intermediate value theorem gives
an exact root for every sufficiently large even $N$.
Strict monotonicity of $g_1$ identifies the unique limiting aspect
ratio, and the first equation determines the positive limit of $b_N$.
\end{proof}

\subsection{Positivity of the exact polygon configuration}

\begin{lemma}\label{geo:cell-limit}
Choose the roots from Lemmas~\ref{geo:noncritical-root},
\ref{geo:critical-root}, and \ref{geo:subcritical-root}, and let
$\alpha_*$ denote the corresponding limit. In logarithmic polar
coordinates $r=e^{h\xi}$ and $\theta=h\eta$,
\begin{equation}\label{geo:cell-convergence}
 h^{s-1}F(r\cos\theta,r\sin\theta,0)
 \longrightarrow\frac1{2\pi}D_s^{\alpha_*}(\xi,\eta)
\end{equation}
locally uniformly off the lattice. The convergence is uniform on a
closed centered fundamental cell after deleting a fixed neighborhood
of its source. Near a source, in Cartesian coordinates on the plane,
\begin{equation}\label{geo:near-source}
 h^{s-1}F(v_0+hY,0)=\frac1{2\pi}|Y|^{-s}+G_N(Y),
 \qquad |G_N(Y)|\le C|Y|^2
\end{equation}
for fixed sufficiently small $|Y|$, with $C$ independent of large $N$.
\end{lemma}

\begin{proof}
Put $X_0=(v_0,0)$, $q_j=t^j$, and
\[
 \Psi_j(X)=q_j^{\beta_N}
 [P_N(X/q_j)-1-q_j^{-2}H_2(X)].
\]
In $\Psi_0$ remove the entire self-pole $N^{-1}|X-X_0|^{-s}$;
the other shells are unchanged. Denote these regularized shells by
$\Psi_j^{\rm reg}$. For each fixed $N$, their sum and its derivatives
converge near $X_0$; this neighborhood may depend on $N$.
For a regular function $f$, write
\[
 \mathcal R_hf(Y)=f(X_0+hY)-f(X_0)-h\nabla f(X_0)\cdot Y.
\]
Both matching conditions give the exact identity
\begin{equation}\label{geo:cell-cancel}
 h^{s-1}F(X_0+hY)=\frac1{2\pi}|Y|^{-s}
          +h^{s-1}\sum_j\mathcal R_h\Psi_j^{\rm reg}(Y).
\end{equation}
The first-order subtraction retains the cancellation needed at $s=1$.

Use centered angular representatives $-N/2\le\ell<N/2$ and set
$v=(-\alpha_Nj,\ell)\in\Lambda_{\alpha_N}$. The scaled source is
\begin{equation}\label{geo:rev-local-map}
 \Phi_h(v)=h^{-1}(e^{hv_1}\cos(hv_2)-1,e^{hv_1}\sin(hv_2)).
\end{equation}
Fix a sufficiently small $d_0>0$. For local indices $h|v|\le d_0$,
$|\Phi_h(v)-v|\le Ch|v|^2$ and $|\Phi_h(v)|\asymp|v|$,
uniformly as $\alpha_N$ stays in a compact positive interval.
Let $K$ be a compact set of observation points disjoint from the
limiting lattice, and put
\[
 \mathcal Q(Y,w)=|Y-w|^{-s}-|w|^{-s}-\nabla g(-w)\cdot Y,
 \qquad g(w)=|w|^{-s}.
\]
For large $|w|$ and $Y\in K$, Taylor's formula and the kernel
derivative bounds yield
\[
 |\mathcal Q(Y,w)|\le C_K|w|^{-s-2},\qquad
 |\nabla_w\mathcal Q(Y,w)|\le C_K|w|^{-s-3}.
\]
Consequently
\begin{equation}\label{geo:rev-kernel-difference}
 |\mathcal Q(Y,\Phi_h(v))-\mathcal Q(Y,v)|
                    \le C_Kh|v|^{-s-1}.
\end{equation}
The normalized point contribution to \eqref{geo:cell-cancel} is
$(2\pi)^{-1}q_j^{s+\beta_N}\mathcal Q(Y,\Phi_h(v))$.
Since $q_j^{s+\beta_N}=1+O(h|v|)$ on the local range, the weight
error has the same bound. The total local error is thus at most
\begin{equation}\label{geo:rev-local-error}
 C_Kh\sum_{1\le|v|\le d_0/h}|v|^{-s-1}
 \le\begin{cases}
 C_Kh\log(2/h),&s=1,\\
 C_Kh,&s>1.
 \end{cases}
\end{equation}
The finitely many remaining small indices converge directly, and the
tail of the limiting series beyond $|v|=L$ is $O_K(L^{-s})$.

For distant shells retain the complete subtraction in $\Psi_j$.
The even Taylor remainder at the origin and the direct kernel bound
give, uniformly near $X_0$,
\begin{equation}\label{geo:rev-distant-brackets}
 |D_X^2\Psi_j(X)|\le
 \begin{cases}
 Cq_j^{\beta_N-4},&q_j\ge2,\\
 C(q_j^{s+\beta_N}+q_j^{\beta_N-2}),&q_j\le1/2.
 \end{cases}
\end{equation}
Since the derivative normalization is $h^{s+1}$ and $T\asymp h$,
the total first-order remainder of these shells is bounded by
\begin{equation}\label{geo:rev-distant-sums}
 C_Kh^s\left(1+\frac1{\beta_N-2}+\frac1{4-\beta_N}\right).
\end{equation}
For $1/2<q_j<2$, nonlocal sources have a fixed positive distance
from $X_0$. Their averaged Hessians are bounded per shell, and there
are $O(h^{-1})$ shells. Their contribution, together with the local
quadratic counterterms, is $O_K(h^s)$. Outer quadratic terms have
not been summed separately.

For all three root sequences,
\[
 \frac{h^s}{\beta_N-2}+\frac{h^s}{4-\beta_N}\longrightarrow0.
\]
Indeed, $\beta_N\to3$ when $s=1$,
$\beta_N-2\asymp1/\log(1/h)$ when $s=2$, and
$\beta_N-2\asymp h^{s-2}$ when $s>2$.
Equations \eqref{geo:rev-local-error}--\eqref{geo:rev-distant-sums}
prove the Cartesian version of \eqref{geo:cell-convergence}.
In logarithmic coordinates $X=X_0+hY_h$ with
$Y_h\to(\xi,\eta)$ uniformly on compact sets, giving the stated form.

After removing the exact pole, the same Hessian estimates give
\begin{equation}\label{geo:rev-uniform-hessian}
 \|D_Y^2G_N\|_{L^\infty(B_\varepsilon)}
 \le C\left[\sum_{v\in\Lambda_{\alpha_N}\setminus\{0\}}|v|^{-s-2}
    +h^s+\frac{h^s}{\beta_N-2}+\frac{h^s}{4-\beta_N}\right]\le C,
\end{equation}
where $\varepsilon>0$ is smaller than one quarter of the uniform
lattice separation. The matching equations imply
$G_N(0)=\nabla G_N(0)=0$, proving \eqref{geo:near-source}.
Finally, the centered cell $|\xi|\le\alpha_N/2$, $|\eta|\le1/2$
contains only its central source and is uniformly separated from
neighboring sources. The compact-set convergence applies to the
entire cell outside a fixed small neighborhood of that source.
\end{proof}

\begin{proof}[Proof of Proposition~\ref{geo:seed}]
Choose the exact roots given above. The limiting aspect ratio lies in
the interval of strict positivity in Lemma~\ref{geo:lattice-positive}.
On a fundamental cell away from its source neighborhoods,
Lemma~\ref{geo:cell-limit} gives a strict positive lower bound.
In those neighborhoods, the positive singular term in
\eqref{geo:near-source} dominates $C|Y|^2$ after fixing the neighborhood
small enough. The logarithmic cell covers one radial period and one
angular period. Since
\[
 F(tX)=t^\beta F(X),
\]
we obtain \eqref{geo:floor} on the whole punctured plane.
The same scaling, polygon symmetry, and $O(n-2)$ symmetry propagate
the two matching equations to the value and full gradient at every source.

It remains to extend positivity off the plane. On $t\le|X|\le1$,
$F$ has a positive lower bound near the source plane, including near
the poles, where it tends to $+\infty$. The negative part of $F$ is
therefore supported away from $z=0$, and $(-F)_+/|z|^2$ is bounded.
Choose $\kappa$ larger than this bound and increase it once more to
obtain a strict lower bound for $|X|^\beta+|z|^2$ on that annulus.
For inner annuli, the scaling identity and $\beta>2$ preserve this
bound: the positive quadratic term decays more slowly than $F$.
This proves \eqref{geo:ambient-positive}.
\end{proof}

\section{The shell parameters and a right inverse}\label{sym:section}

\subsection{The exact linearization}

Fix a configuration from Proposition~\ref{geo:seed}. We shall increase
$N$ once more below. Write $V=\{v_0,\ldots,v_{N-1}\}\subset S^1$
for its vertices.
Set $r_0=\sqrt t$ and $R_i=r_0t^i$ for $i\in\Z$. For parameters $\vartheta_i=(\eta_i,b_i)$, define
\[
 M_i=R_i^{s+\beta}e^{\eta_i},\qquad r_i=R_ie^{b_i},\qquad
 q_{i,\nu}=r_i(v_\nu,0).
\]
Let $T_2$ denote the Taylor polynomial of degree two at the origin and put $N_2=I-T_2$. The field associated with these parameters is
\begin{equation}\label{sym:field}
 F_\vartheta(X)=\sum_{i\in\Z}\mathcal P_i(X),\qquad
 \mathcal P_i(X)=\frac{M_i}{N}\sum_{\nu=0}^{N-1}N_2|X-q_{i,\nu}|^{-s}.
\end{equation}
The terms of each shell are combined before summation. If the parameters are uniformly small, the inner and outer tails converge with rates $t^{\beta-2}$ and $t^{4-\beta}$, respectively. At $q_i=q_{i,0}$, let $F_{\vartheta,i}^{\mathrm{reg}}$ be obtained by removing the complete kernel $(M_i/N)|X-q_i|^{-s}$. Define
\begin{equation}\label{sym:map}
 \mathcal F_i(\vartheta)=
 \begin{pmatrix}
 R_i^{-\beta}F_{\vartheta,i}^{\mathrm{reg}}(q_i)\\
 R_i^{-\beta}q_i\cdot\nabla F_{\vartheta,i}^{\mathrm{reg}}(q_i)
 \end{pmatrix}.
\end{equation}
The other components of the gradient vanish by symmetry. The corresponding balance equations give $\mathcal F(0)=0$.

For the linearization, introduce
\[
 K(r,q)=\frac1N\sum_{\nu=0}^{N-1}
       (r^2+q^2-2rq\,v_0\cdot v_\nu)^{-s/2},\qquad
 c=\frac{s^2}{4}.
\]
When the source and observation belong to the same shell, the term $\nu=0$ is deleted from $K$ and all its derivatives. For $q=t^{j-i}$ set
\begin{align}
 A(q)&=q^{s+\beta}K(1,q)-q^\beta-cq^{\beta-2},\notag\\
 B(q)&=q^{s+\beta+1}K_q(1,q)+sq^\beta+(s+2)cq^{\beta-2},\notag\\
 C(q)&=q^{s+\beta}K_r(1,q)-2cq^{\beta-2},\label{sym:blocks}\\
 D(q)&=q^{s+\beta+1}K_{rq}(1,q)+2c(s+2)q^{\beta-2},\notag\\
 E(q)&=q^{s+\beta}K_{rr}(1,q)-2cq^{\beta-2}.\notag
\end{align}
Direct differentiation of \eqref{sym:map} gives
\begin{align}
 (D\mathcal F(0)\vartheta)_{i,1}
 &=\sum_j\bigl[A(t^{j-i})\eta_j+B(t^{j-i})b_j\bigr],\notag\\
 (D\mathcal F(0)\vartheta)_{i,2}
 &=\sum_j\bigl[C(t^{j-i})\eta_j+D(t^{j-i})b_j\bigr]
     +b_i\sum_jE(t^{j-i}).\label{sym:linear}
\end{align}
The observation point contributes no additional term to the first row, since its regular radial derivative vanishes at the base configuration. In the second row, differentiating $r\partial_rF$ gives $r\partial_rF+r^2\partial_{rr}F$; the first term vanishes after summation, and the second is the last term of \eqref{sym:linear}.

Write this real convolution operator as
\[
 (\mathcal A\vartheta)_i=\sum_{j\in\Z}A_{i-j}\vartheta_j,
 \qquad \mathbf A(z)=\sum_{k\in\Z}A_kz^k.
\]
The matrices satisfy
\begin{equation}\label{sym:decay}
 \|A_k\|\le C
 \begin{cases}
 t^{(4-\beta)k},&k\ge1,\\
 t^{(\beta-2)(-k)},&k\le-1.
 \end{cases}
\end{equation}
Indeed, for $q>1$ the constant and quadratic Taylor terms cancel, leaving $O(q^{\beta-4})$, together with all parameter derivatives used in \eqref{sym:blocks}. For $q<1$, the slowest term is $q^{\beta-2}$. Thus $\mathbf A$ is holomorphic in
\begin{equation}\label{sym:annulus}
 t^{\beta-2}<|z|<t^{-(4-\beta)}.
\end{equation}

\subsection{A nonzero value of the determinant}

We prove $\det\mathbf A(-1)>0$ for sufficiently large $N$ in every
dimension. This single nonzero value is enough for the right inverse.

Use the coefficients $C_d=C_{d,N}$ of Lemma~\ref{geo:coeff}; all sums indexed by $d$ below run over even integers. Set
\[
 S_N=\frac1N\sum_{\nu=1}^{N-1}(2\sin(\pi\nu/N))^{-s},\quad
 \delta=\beta-2,\quad u_d=s+\beta+d,\quad v_d=d-\beta,
\]
and define
\[
 H_u=\frac1{e^{Tu}-1},\qquad J_u=-\frac1{e^{Tu}+1},\qquad T=\alpha h.
\]
Write $a_{ij}=(\mathbf A(-1))_{ij}$. Summing the shell index first in \eqref{sym:blocks} yields
\begin{align}
 a_{11}={}&S_N-1-c+\sum_{d\ge0}C_dJ_{u_d}
                 -J_\beta-cJ_\delta+\sum_{d\ge4}C_dJ_{v_d},
 \label{sym:pi11}\\
 a_{12}={}&-\frac{s}{2}S_N+s+(s+2)c
       +\sum_{d\ge0}dC_dJ_{u_d}+sJ_\beta+(s+2)cJ_\delta\notag\\
       &\hspace{22mm}-\sum_{d\ge4}(s+d)C_dJ_{v_d},\label{sym:pi12}\\
 a_{21}={}&-\frac{s}{2}S_N-2c
       -\sum_{d\ge0}(s+d)C_dJ_{u_d}-2cJ_\delta
       +\sum_{d\ge4}dC_dJ_{v_d},\label{sym:pi21}\\
 a_{22}={}&\frac{s(s+1)}2S_N+2c(s+1)
       -2cH_\delta+2c(s+2)J_\delta\notag\\
       &+\sum_{d\ge0}C_d\bigl[(s+d)(s+d+1)H_{u_d}
                                      -d(s+d)J_{u_d}\bigr]\notag\\
       &+\sum_{d\ge4}C_d\bigl[d(d-1)H_{v_d}
                                      -d(s+d)J_{v_d}\bigr].\label{sym:pi22}
\end{align}
These series are absolutely convergent for each fixed $N,T,\beta$. To check the same-shell term in the last formula, symmetry and homogeneity give
\[
 K_r(1,1)=-\frac{s}{2}S_N,
 \qquad K_{rr}(1,1)+K_{rq}(1,1)=\frac{s(s+1)}2S_N.
\]
The observation and source derivatives must therefore be added before taking the limit. Their separate highest-order terms cancel. The corresponding quadratic counterterms add to $2c(s+1)$.

Homogeneity also gives $B(q)+C(q)=-sA(q)$. Consequently,
\begin{equation}\label{sym:detidentity}
 a_{12}+a_{21}=-sa_{11},\qquad
 \det\mathbf A(-1)
 =a_{11}\left(a_{22}-\frac{s^2}{4}a_{11}\right)
       +\frac{(a_{12}-a_{21})^2}{4}.
\end{equation}
We shall use the positive last term without estimating its limit.

\begin{lemma}\label{sym:nonzero}
For every fixed integer $s\ge1$, the balanced polygons of Proposition~\ref{geo:seed} satisfy $\det\mathbf A(-1)>0$ when $N$ is sufficiently large.
\end{lemma}
\begin{proof}
First suppose $s\ge2$. Let $\alpha_N\to\alpha$ and $\beta_N\to2$ be the parameters constructed there. Define
\[
 G_s(y)=\sum_{\ell\in\Z}(y^2+\ell^2)^{-s/2},\quad G_0=2\zeta(s),
\]
and put
\[
 \widehat G_s(\pi)=G_0+2\sum_{j\ge1}(-1)^jG_s(\alpha j),\qquad
 L_s(\pi)=4\sum_{\substack{j\ge1\\j\ \mathrm{odd}}}G_s''(\alpha j).
\]
The first series is understood by its ordinary alternating sum when $s=2$.

By Lemma~\ref{geo:coeff}, if $d/N\to x>0$ away from the integers, then
$N^{2-s}C_{d,N}\to W_s(x)$, and
\[
 N^{2-s}C_{d,N}\le C_s(x+N^{-1})^{s-2}(1+x).
\]
The even $d$ have mesh $2/N$. Applying this bound to \eqref{sym:pi11}, with the finitely many low terms separated, gives
\begin{align}
 h^{s-1}a_{11}
 &\longrightarrow\frac{G_0}{2\pi}
 -(2\pi)^{s-1}\int_0^\infty
        \frac{W_s(x)}{e^{2\pi\alpha x}+1}\dd x
 =\frac{\widehat G_s(\pi)}{2\pi}.\label{sym:limit11}
\end{align}
Near zero the integrand is $O(x^{s-2})$, and at infinity it has exponential decay. The discrete bound controls these two tails uniformly, so the Riemann limit holds also for $s=2$. The last equality follows from the Laplace formula for $G_s$ in Lemma~\ref{geo:coeff}; an Abel factor justifies summing the alternating series under the integral at the endpoint $s=2$.

Since $H_u-J_u=1/\sinh(Tu)$, the two positive principal sums in \eqref{sym:pi22} give
\begin{equation}\label{sym:limit22}
 h^{s+1}a_{22}\longrightarrow
 (2\pi)^{s+1}\int_0^\infty
       \frac{x^2W_s(x)}{\sinh(2\pi\alpha x)}\dd x
 =\frac{L_s(\pi)}{2\pi}>0.
\end{equation}
Here the integrand is $O(x^{s-1})$ at zero. The same-shell term contributes $O(h^2)$. The only possibly large negative term satisfies
\[
 h^{s+1}H_\delta\le C\frac{h^s}{\delta}\longrightarrow0:
\]
for $s>2$ one has $\delta\asymp h^{s-2}$, while for $s=2$ one has $\delta\asymp1/\log N$. All remaining bounded or lower-order terms vanish. The positivity in \eqref{sym:limit22} also follows directly from the positive Laplace density $W_s$.

It remains to show $\widehat G_s(\pi)>0$. Suppose first that $s\ge3$. The Poisson expansion, with $\nu=(s-1)/2$, is
\[
 G_s(y)=\frac{\sqrt\pi\Gamma(\nu)}{\Gamma(s/2)}y^{1-s}
 +\frac{4\pi^{s/2}}{\Gamma(s/2)}y^{-\nu}
       \sum_{\ell\ge1}\ell^\nu K_\nu(2\pi\ell y).
\]
The identity $(y^\nu K_\nu(cy))'=-cy^\nu K_{\nu-1}(cy)$ implies
$(y^{s-1}G_s(y))'<0$. At the pressure root of Lemma~\ref{geo:pressure},
\[
 0=G_0+2\sum_{j\ge1}\bigl[G_s(\alpha j)+\alpha jG_s'(\alpha j)\bigr].
\]
It follows that
\[
 G_0>2(s-2)\sum_{j\ge1}G_s(\alpha j)
       \ge2\sum_{j\ge1}G_s(\alpha j),
\]
and hence $\widehat G_s(\pi)>0$.

For $s=2$, Lemma~\ref{geo:critical-root} gives $6/\pi<\alpha<2$. Using $G_2(y)=\pi\coth(\pi y)/y$, we obtain
\begin{align*}
 \widehat G_2(\pi)
 &=\frac{\pi^2}{3}-\frac{2\pi}{\alpha}\log2
   -\frac{4\pi}{\alpha}\sum_{k\ge1}\log(1+e^{-2\pi\alpha k})\\
 &\ge\frac{\pi^2}{3}(1-\log2)
      -\frac{2\pi^2}{3}\frac{e^{-12}}{1-e^{-12}}>0.
\end{align*}
Thus \eqref{sym:limit11}--\eqref{sym:limit22} imply
$a_{11}>0$ and $a_{22}-(s^2/4)a_{11}>0$ for all sufficiently large $N$. The conclusion follows from \eqref{sym:detidentity}.

It remains to consider $s=1$. By Lemma~\ref{geo:subcritical-root}, for
all large $N$ the parameters lie in
\[
 5/2\le\alpha\le27/10,\qquad 3\le\beta\le7/2.
\]
Use $K_N,K,\delta_N,f_\beta$ and $\sigma=\beta+1/2$ from its proof.
The alternating version of \eqref{geo:s1-pair-value} is
\begin{equation}\label{sym:s1-alternating}
 \begin{aligned}
 a_{11}={}&S_N-\frac54+\sum_{j\ge1}(-1)^j f_\beta(jT)\\
       &+2\sum_{j\ge1}(-1)^j\cosh(\sigma jT)\delta_N(jT).
 \end{aligned}
\end{equation}
For the angular average, \eqref{geo:log-euler} gives
\[
 \sum_{j\ge1}f_\beta(jT)
 =\frac1T\int_0^\infty f_\beta(u)\dd u
       +\frac1\pi\log\frac{T}{16\pi}+\frac54+o(1).
\]
Since the alternating sum is twice the sum with step $2T$ minus the
sum with step $T$, its integral term cancels and its remaining terms
are $\pi^{-1}\log(T/(4\pi))+5/4+o(1)$. The cosecant asymptotic in
Lemma~\ref{geo:subcritical-root} therefore supplies
$(\gamma+\log(\alpha/\pi))/\pi$ after adding $S_N-5/4$.
For the nonzero angular modes,
\eqref{geo:s1-fourier-limit}--\eqref{geo:s1-fourier-bound} apply
unchanged: multiplication by $(-1)^j$ preserves absolute domination.
Thus, uniformly on the displayed parameter rectangle,
\begin{equation}\label{sym:s1-limit11}
 \begin{aligned}
 a_{11}&=f_1(\alpha)+o(1),\\
 f_1(\alpha)&=\frac{\gamma+\log(\alpha/\pi)}\pi
       +\frac4\pi\sum_{j,k\ge1}(-1)^jK_0(2\pi\alpha jk).
 \end{aligned}
\end{equation}
This argument does not use the coefficient bound for $W_s$, which
was restricted to $s\ge2$.

The integral representation of $K_0$, with $\cosh v\ge1+v^2/2$, gives
$K_0(u)\le\sqrt{\pi/(2u)}e^{-u}$. Set $q=e^{-2\pi\alpha}$.
The absolute value of the series correction in \eqref{sym:s1-limit11}
is at most
\[
 \frac2{\pi\sqrt\alpha}\sum_{j,k\ge1}q^{jk}
 \le\frac2{\pi\sqrt\alpha}\frac{q}{(1-q)^2}<10^{-6},
\]
where $jk\ge j+k-1$ was used. Since $\gamma>0.57$ and
$\log(\alpha/\pi)>-0.23$ on this interval, the first term of
$f_1$ exceeds $0.34/\pi>0.108$. Hence $f_1>1/10$.

For the second diagonal entry, all terms in both sums of
\eqref{sym:pi22} are nonnegative. If $d$ is even and
$N/4\le d\le N/2$, then
\[
 C_{d,N}=b_{d/2}^2\ge\frac1{2d},\qquad b_m=(1/2)_m/m!.
\]
The inequality follows by induction from $b_m^2\ge1/(4m)$ for
$m\ge1$. For $N\ge100$, there are at least $N/10$ such indices,
and $T(1+\beta+d)<10$. Retaining only these terms in the first sum
of \eqref{sym:pi22} gives a contribution of at least
$N^2/[80(e^{10}-1)]$. The remaining terms outside the sums satisfy
\[
 -\tfrac12H_\delta\ge-\frac1{2T(\beta-2)}
       \ge-\frac{N}{10\pi},\qquad
 S_N+1+\tfrac32J_\delta\ge\frac14.
\]
Consequently
\begin{equation}\label{sym:s1-bound22}
 a_{22}\ge\frac{N^2}{80(e^{10}-1)}-\frac{N}{10\pi}.
\end{equation}
Equations \eqref{sym:s1-limit11}--\eqref{sym:s1-bound22} give
$a_{11}>0$ and $a_{22}-a_{11}/4>0$ for all large $N$.
Applying \eqref{sym:detidentity} proves the assertion also for $s=1$.

\end{proof}

Fix one such $N$, together with $t$ and $\beta$. All subsequent
constants may depend on this fixed configuration.

\subsection{A weighted inverse with finitely many exterior parameters}

For $\rho>1$, write
\[
 \mathcal E_\rho=\left\{x=(x_i)_{i\in\Z}:\quad
        \|x\|_\rho=\sum_{i\in\Z}\rho^i|x_i|<\infty\right\}.
\]
Let $P$ retain the indices $i\ge0$, let $P_M$ retain $i\ge-M$, and put $Q_M=I-P_M$. Set $Y_\rho=P\mathcal E_\rho$ and $X_{\rho,M}=P_M\mathcal E_\rho$. All these spaces consist of real two-component sequences.

\begin{proposition}\label{sym:rightinverse}
The number $\rho>1$ can be chosen arbitrarily close to $1$ so that, for some finite $M$, the operator
\[
 A_+=P\mathcal A\big|_{X_{\rho,M}}:X_{\rho,M}\longrightarrow Y_\rho
\]
has a bounded right inverse $\mathcal R$. Its range is a closed complemented subspace of $X_{\rho,M}$.
\end{proposition}
\begin{proof}
By Lemma~\ref{sym:nonzero}, the holomorphic function $\det\mathbf A$ is not identically zero in \eqref{sym:annulus}. Its zeros are isolated. We may therefore choose $\rho>1$, as close to $1$ as needed, such that $\det\mathbf A$ has no zero on $|z|=\rho$. In particular, choose $\rho<t^{-(4-\beta)}$. The inverse matrix function is holomorphic on a smaller annulus containing this circle. Write
\[
 \mathbf A(z)^{-1}=\sum_{k\in\Z}B_kz^k.
\]
Cauchy's estimates on two circles of radii smaller and larger than $\rho$ give
\begin{equation}\label{sym:wiener}
 \sum_k\rho^k\|A_k\|<\infty,\qquad
 \sum_k\rho^k\|B_k\|<\infty.
\end{equation}
The convolution $\mathcal B$ with coefficients $B_k$ is consequently a two-sided inverse of $\mathcal A$ on $\mathcal E_\rho$. These statements follow from absolute multiplication of the Laurent series and the weighted Young inequality. The coefficients are real, by complex conjugation.

For $f\in Y_\rho$, extend $f$ by zero to negative indices and define
\[
 \mathcal R_M=P_M\mathcal BP,\qquad
 E_M=P\mathcal A Q_M\mathcal BP.
\]
Then $P\mathcal A\mathcal R_M=I-E_M$. Moreover,
\begin{align*}
 \|Q_M\mathcal BPf\|_\rho
 &\le\sum_{j\ge0}\rho^j|f_j|
                  \sum_{k<-M-j}\rho^k\|B_k\|\\
 &\le\left(\sum_{k<-M}\rho^k\|B_k\|\right)\|f\|_\rho.
\end{align*}
It follows that
\[
 \|E_M\|\le
 \left(\sum_k\rho^k\|A_k\|\right)
 \left(\sum_{k<-M}\rho^k\|B_k\|\right)\longrightarrow0.
\]
Choose a finite $M$ with $\|E_M\|<1$ and set
\begin{equation}\label{sym:Rformula}
 \mathcal R=P_M\mathcal BP(I-E_M)^{-1}.
\end{equation}
This is a bounded right inverse. Since $A_+\mathcal R=I$, the operator $\mathcal R A_+$ is a bounded projection onto $\operatorname{Ran}\mathcal R$. The range is therefore closed and complemented.
\end{proof}

The direction of the weight is relevant: conjugation by $x_i\mapsto\rho^ix_i$ replaces the convolution coefficients by $\rho^kA_k$ and the symbol by $\mathbf A(\rho z)$. Thus $\rho>1$ imposes decay toward the inner shells. It may also be chosen to satisfy any fixed finite list of conditions $\rho t^\gamma<1$ with $\gamma>0$. No bound uniform as $N\to\infty$ is needed.

\subsection{The nonlinear point-source map}

In \eqref{sym:field}, keep the parameters with $i<-M$ equal to zero. The corresponding base sources are retained. Only the masses and radii of the shells with $i\ge-M$ are varied. Since $R_{-1}=t^{-1/2}>1$, the $M$ additional shells are outside $B_1$; small parameter changes keep them there. Their fields are smooth and biharmonic in $B_1$.

\begin{proposition}\label{sym:regularity}
The map $\mathcal F_+=(\mathcal F_i)_{i\ge0}$ is $C^2$ from a neighborhood of zero in $X_{\rho,M}$ to $Y_\rho$, and $D\mathcal F_+(0)=A_+$. If
\[
 \mathcal Q(\vartheta)=\mathcal F_+(\vartheta)-A_+\vartheta,
\]
then, on a sufficiently small fixed ball,
\begin{align}
 \|\mathcal Q(\vartheta)\|_\rho&\le C\|\vartheta\|_\rho^2,\notag\\
 \|\mathcal Q(\vartheta)-\mathcal Q(\vartheta')\|_\rho
 &\le C(\|\vartheta\|_\rho+\|\vartheta'\|_\rho)
                          \|\vartheta-\vartheta'\|_\rho.
 \label{sym:quadratic}
\end{align}
\end{proposition}
\begin{proof}
The embedding $\|\vartheta\|_\infty\le\rho^M\|\vartheta\|_\rho$ permits us to impose $\sup_i|b_i|<T/4$. Different shells then remain separated. After normalization by $R_i^\beta$, write row $i$ of \eqref{sym:map} as
\[
 \mathcal F_i(\vartheta)=\sum_{j\in\Z}
            \mathcal K_{j-i}(\vartheta_i,\vartheta_j),
\]
where the same-shell summand uses the deleted kernel. The first three derivatives of each summand, with respect to its two vector arguments, satisfy
\begin{equation}\label{sym:C3kernel}
 \|D^\ell\mathcal K_k\|\le C d_k\quad(\ell=1,2,3),\qquad
 d_k=\begin{cases}
 t^{(\beta-2)k},&k\ge1,\\
 1,&k=0,\\
 t^{(4-\beta)(-k)},&k\le-1.
 \end{cases}
\end{equation}
For outer sources this is the Taylor cancellation through degree two; evenness cancels degree three. Inner sources are controlled by their quadratic counterterm. Parameter differentiation adds bounded exponential factors and preserves these rates. The same-shell deleted kernel is smooth on the fixed parameter neighborhood. In particular,
\[
 \sum_kd_k<\infty,\qquad \sum_kd_k\rho^{-k}<\infty.
\]
For example, the second derivative obeys
\[
 |D^2\mathcal F_i(\vartheta)[v,w]|
 \le C\sum_jd_{j-i}(|v_i|+|v_j|)(|w_i|+|w_j|).
\]
Bound one factor by $2\rho^M\|v\|_\rho$ and sum the other over $i\ge0$ with weight $\rho^i$. The two displayed summability conditions give
\[
 \|D^2\mathcal F_+(\vartheta)[v,w]\|_\rho
       \le C_M\|v\|_\rho\|w\|_\rho.
\]
The first derivative is treated in the same way, and the third derivative bound gives continuity of the second derivative. Since $\mathcal F_+(0)=0$, integration along line segments proves that the map takes values in $Y_\rho$ and is $C^2$. Taylor's formula gives \eqref{sym:quadratic}.
\end{proof}

The exterior variations start at degree four in $X$. Their contribution is therefore $O(\|\vartheta\|_\rho|X|^4)$, with a constant depending on the fixed $M$. In particular, they preserve the positive bound on the source subspace of order $|X|^\beta$ when the parameter ball is sufficiently small.

\section{The approximate solution}\label{bg:section}

\subsection{Renormalized potentials and moving coordinates}
The fundamental solution of $\Delta^2$ in $\R^n$ is
\[
\Gamma(X)=\frac{|X|^{-s}}{2s(s+2)|S^{n-1}|}.
\]
For a function smooth near the origin, let $T_2$ denote its Taylor
polynomial of degree two and put $N_2=I-T_2$. Define
\begin{equation}\label{bg:S}
\mathcal Sf(X)=\int_{B_1}K_2(X,Y)f(Y)\dd Y,\qquad
K_2(X,Y)=N_{2,X}\Gamma(X-Y).
\end{equation}
Here $N_{2,X}$ acts at $X=0$ with $Y\ne0$ fixed. If $f$ is locally
bounded off zero and $\int_{B_1}|Y|^{-s-2}|f(Y)|\dd Y<\infty$, the integral
converges for $X\ne0$, is locally integrable, and satisfies
$\Delta^2\mathcal Sf=f$ distributionally. This follows from local
integrability of $\Gamma$, absolute integrability of the Taylor
coefficients, and Fubini's theorem. Sources on $B_1$ are extended by
zero in convolutions; full-space integrals are specified.

For even sources, use
$K_{2,e}(X,Y)=[K_2(X,Y)+K_2(X,-Y)]/2$. Taylor's formula gives
\begin{equation}\label{bg:kernel}
|D_X^jK_{2,e}(X,Y)|\le C
\begin{cases}
|X|^{4-j}|Y|^{-s-4},& |X|\le |Y|/2,\\
|X|^{2-j}|Y|^{-s-2},& |Y|\le |X|/2,
\end{cases}
\end{equation}
for $0\le j\le2$. At comparable scales retain also the integrable
singularities $|X-Y|^{-s-j}$ and $|X+Y|^{-s-j}$.
Far-shell estimates are applied to complete even shells.

Fix the configuration of Proposition~\ref{geo:seed} and the right inverse of
Proposition~\ref{sym:rightinverse}. We write $X=X_{\rho,M}$ for the parameter
space and $\mathcal F=\mathcal F_+$ for its interior point-source map.
Take $r_0=\sqrt t$, so that
$R_i=r_0t^i$ and the unperturbed shells $i<0$ lie outside $B_1$.
Let $\mathcal P_i$ be the renormalized shell potential from
Sect.~\ref{sym:section}, and put $F_\vartheta=\sum_{i\in\Z}\mathcal P_i$.
Parameters with $i<-M$ remain zero; their original sources are retained.
Use $\mathfrak G=D_N\times O(n-2)$, where $D_N$ is the dihedral
group of the polygon. Since $N$ is even, $-I\in\mathfrak G$.
The stabilizer of a
vertex fixes only its radial direction, so only scale and radial
translation Jacobi fields remain in the invariant class.

Set
\begin{equation}\label{bg:Q}
A_n=[n(n-4)(n^2-4)]^{s/8},\qquad Q(y)=A_n(1+|y|^2)^{-a}.
\end{equation}
Then $\Delta^2Q=Q^p$. For $i\ge0$, define
\begin{equation}\label{bg:mu}
A_n\mu_i^a=\frac{\lambda M_i}{N},\qquad
Q_{i,\nu}(X)=\mu_i^{-a}Q\left(\frac{X-q_{i,\nu}}{\mu_i}\right).
\end{equation}
Thus $\mu_i\asymp\lambda^{2/s}R_i^{2+2\beta/s}$ on a fixed small
parameter ball. Put $\Psi_i=\sum_\nu Q_{i,\nu}$ and
\begin{equation}\label{bg:W}
W_\vartheta=\lambda\kappa|z|^2+
 \lambda\sum_{i<0}\mathcal P_i+\sum_{i\ge0}N_2\Psi_i.
\end{equation}
The constant $\kappa$ is fixed below. The series converges locally
smoothly in $B_1\setminus\{0\}$, including at each center, and
\begin{equation}\label{bg:WPDE}
\Delta^2W_\vartheta=\sum_{i\ge0,\nu}Q_{i,\nu}^p.
\end{equation}

When convenient, an index $i$ also includes a vertex, with its shell
radius still denoted by $R_i$. Choose $\delta>0$ so that the reference
balls $B(q_i(0),20\delta R_i)$ are pairwise disjoint and contained in
$B_1$. Restrict the parameter ball so that the moving caps remain
separated. Define
\begin{equation}\label{bg:T}
T_i^\vartheta(X)=q_i(\vartheta)+\sigma_i(\vartheta)(X-q_i(0)),
\qquad \sigma_i=\frac{\mu_i(\vartheta)}{\mu_i(0)}=e^{\eta_i/a},
\end{equation}
and $\Omega_i(j)=T_i^\vartheta(B(q_i(0),j\delta R_i))$.
Take smooth radial cutoffs $\chi,\zeta$ with values in $[0,1]$:
$\chi=1$ on $B_{1/2}$, $\supp\chi\subset B_1$,
$\zeta=1$ on $B_2$, and $\supp\zeta\subset B_3$.
Define $\chi_i,\zeta_i$ by evaluating these functions at
$((T_i^\vartheta)^{-1}X-q_i(0))/(\delta R_i)$.
The micro-radius $L_i=\delta R_i/\mu_i(0)$ is independent of the
parameters and tends uniformly to infinity as $\lambda\to0$.

Using fixed cutoffs on the reference caps, extend $T_i^\vartheta$ to a
$\mathfrak G$-equivariant map $\Phi_\vartheta:B_1\to B_1$, equal to
$T_i^\vartheta$ on the $10\delta R_i$ balls and to the identity outside
the $20\delta R_i$ balls. It is a bi-Lipschitz diffeomorphism off zero
and fixes zero continuously. For $\epsilon=\|\vartheta-\vartheta'\|_X$,
\begin{equation}\label{bg:Phi}
|\delta\Phi(X)|\le C\epsilon|X|,\quad
\|D\delta\Phi\|_\infty\le C\epsilon,\quad
|D^2\delta\Phi(X)|\le C\epsilon/|X|.
\end{equation}
Indeed, the displacement is $O(\epsilon R_i)$ on each reference cap,
and each derivative of the cutoff costs $R_i^{-1}$.
Taking $\|D\Phi-I\|_\infty<1/2$ gives invertibility.

\subsection{Positivity and fine background estimates}
Use annuli $\mathcal A_0=\{t<|X|<1\}$ and
$\mathcal A_i=\{R_it^{1/2}<|X|\le R_it^{-1/2}\}$ for $i\ge1$.
Decrease $\delta$ so that the reference $20\delta R_i$ balls and all
moving caps remain in their respective annuli. Then $\Phi_\vartheta$
preserves every annulus. Define
\begin{equation}\label{bg:space}
\|v\|_{\mathcal Y}=\sum_{i\ge0}\rho^iR_i^{-\beta}
       \|v\|_{L^\infty(\mathcal A_i)},\qquad
\mathcal Y=\{v\in C(B_1\setminus\{0\})^{\mathfrak G}:
                                  \|v\|_{\mathcal Y}<\infty\}.
\end{equation}
Then $|v(X)|\le C\|v\|_{\mathcal Y}|X|^\beta$.
For vertex-indexed sequences, take the maximum on each shell before
summation.
Choose the zero-free radius $\rho>1$ sufficiently
close to one that $\rho t^\eta<1$ for every positive shell exponent
$\eta$ used below. There are only finitely many such exponents for
fixed $s$. Pullback by $\Phi_\vartheta$ and its inverse is bounded on
$\mathcal Y$. All these restrictions are imposed when $\rho$ is first
selected in Proposition~\ref{sym:rightinverse}; $M$ and the right inverse
are then fixed, before $\lambda$ is chosen.

\begin{lemma}\label{bg:bounds}
One can first fix $\kappa$, then $\delta$ and a small parameter ball,
and finally $\lambda_0>0$, so that for $0<\lambda<\lambda_0$,
\begin{equation}\label{bg:positive}
W_\vartheta\ge c\lambda|X|^\beta,\qquad
W_\vartheta\ge c(Q_i+\lambda|z|^2)\quad\text{on }\Omega_i(8).
\end{equation}
Outside the $\Omega_i(1/2)$, $|W_\vartheta|\le C\lambda|X|^2$.
For $B_i=W_\vartheta-Q_i$, on $\Omega_i(8)$,
\begin{align}
B_i(q_i)&=\lambda R_i^\beta E_i(\vartheta)
                        +O(\lambda^{1+4/s}R_i^4),\label{bg:value}\\
q_i\cdot\nabla B_i(q_i)&=\lambda R_i^\beta F_i(\vartheta)
                        +O(\lambda^{1+4/s}R_i^4),\label{bg:gradient}\\
|D^jB_i|&\le C\lambda R_i^{2-j}\quad(0\le j\le3).
\label{bg:jets}
\end{align}
Here $(E_i,F_i)$ is the point-source map. The two normalized error
sequences in \eqref{bg:value}--\eqref{bg:gradient} and their parameter
Lipschitz constants in $\ell^1_\rho$ are $O(\lambda^{4/s})$.
Pulled-back derivatives in \eqref{bg:jets}, for $j\le2$, have parameter
differences bounded by the same right-hand side times $\epsilon$.
Moreover,
\begin{equation}\label{bg:fine}
\big|\delta(\sigma_i^aB_i\circ T_i^\vartheta)(X)\big|
 \le C\epsilon\lambda(R_i^\beta+|z|^2).
\end{equation}
In the normalization
$\overline W_i=\sigma_i^aW_\vartheta\circ T_i^\vartheta$ the main profile
is fixed, and $|\delta\overline W_i|\le C\epsilon\overline W_i$.
The analogous relative estimate holds outside the caps.
\end{lemma}

\begin{proof}
On a fundamental annulus the unperturbed field on the source subspace
is positive by Proposition~\ref{geo:seed}, and its poles are positive. Its negative set is
therefore separated from $z=0$. Thus $(-F_0)_+/|z|^2$ is bounded there.
Increasing $\kappa$ and using $\beta>2$ gives
\begin{equation}\label{bg:ambient}
F_0(X)+\kappa|z|^2\ge c(|X|^\beta+|z|^2),\qquad 0<|X|<1.
\end{equation}
The factor $r_0$ changes only fixed constants.
Off the caps, differentiation of the renormalized shell kernels gives
$|D^j(F_\vartheta-F_0)|\le C\|\vartheta\|_X|X|^{\beta-j}$ for $j\le2$.

We detail the regularization error. For an interior point let
$\tau_i=\lambda M_i/N$ and
$\mathcal E_i=Q_i-\tau_i|X-q_i|^{-s}$. For $d=|X-q_i|\ge cR_i$,
the integral remainder for $(1+v)^{-a}$ gives
\begin{equation}\label{bg:soft-tail}
|D^j\mathcal E_i|\le C\tau_i\mu_i^2d^{-s-2-j},\qquad 0\le j\le6.
\end{equation}
Parameter differences in moving coordinates satisfy the same bound
times $\epsilon$. Since
$\tau_i\mu_i^2\asymp\lambda^{1+4/s}R_i^{s+\beta+4+4\beta/s}$,
complete even outer shells, inner shells, and comparable shells
contribute, respectively,
\begin{align}
&C\lambda^{1+4/s}r^{4-j}
       \sum_{r\lesssim R_i\le R_0}R_i^{\beta-2+4\beta/s},\notag\\
&C\lambda^{1+4/s}r^{2-j}
       \sum_{R_i\lesssim r}R_i^{\beta+4\beta/s},\notag\\
&C\lambda^{1+4/s}r^{\beta+2+4\beta/s-j}.
\label{bg:soft-sums}
\end{align}
Only $i\ge0$ is regularized. The first exponent is positive and
$\beta+2+4\beta/s>4$. Hence the total is
$C\lambda^{1+4/s}r^{4-j}$ for $j\le3$ away from the caps.
Inside a cap, delete its own profile and pole. Its own error is then
only $-T_2\mathcal E_i$, satisfying the same bound with $r=R_i$.
Other sources are at a fixed macroscopic distance. This proves the
soft-error bounds, including their parameter differences.

With the self-pole removed, the point field has derivatives bounded
by $CR_i^{\beta-j}$ on the cap. Its central value and radial derivative
are exactly the entries of the point-source map. The normal polynomial
has zero value and gradient there. This proves
\eqref{bg:value}--\eqref{bg:jets}.
The sequence $R_i^{4-\beta}$ belongs to $\ell^1_\rho$, and the point-source
map is $C^2$ in that space. This also proves the weighted central
error and difference bounds.

For the fine difference put
$B_i^{\rm flat}=B_i-\lambda\kappa|z|^2$. Shell differentiation gives
\[
|D^jB_i^{\rm flat}|\le C\lambda R_i^{\beta-j},\qquad
|D^j\delta(B_i^{\rm flat}\circ T_i^\vartheta)|
       \le C\epsilon\lambda R_i^{\beta-j}\quad(j\le2).
\]
The self-profile is absent; no derivative of a concentrated peak
appears. The transported normal polynomial is
$\lambda\kappa\sigma_i^2|z|^2$. Multiplication by $\sigma_i^a$ yields
\eqref{bg:fine}.

Inside the cap, after deleting the nonnegative normal polynomial, the
negative part of the regular background is at most
\[
C\lambda R_i^\beta
 \left(\|\vartheta\|_X+\|\vartheta\|_X\frac d{R_i}
                         +\frac{d^2}{R_i^2}\right)
       +C\lambda^{1+4/s}R_i^4.
\]
Here the first two terms use the vanishing regular value and gradient
at the balanced configuration. Also,
\[
Q_i\ge c\lambda R_i^\beta
 \left[(\mu_i/R_i)^2+(d/R_i)^2\right]^{-s/2}.
\]
Choosing $\delta$, the parameter ball, and then $\lambda$ small makes
the negative part less than $Q_i/2$. Outside the caps apply
\eqref{bg:ambient} and \eqref{bg:soft-sums}. This proves positivity and
the upper bound.
Finally $\sigma_i^aQ_i\circ T_i^\vartheta=Q_i^0$ exactly, and
$Q_i^0\ge c\lambda R_i^\beta$ on the reference cap. The fine difference
and positivity give the relative bound. Off the caps, the separated
shell kernels give
$|\delta W|\le C\epsilon\lambda(r^\beta+|z|^2)$, which is controlled by
the positive lower bound. This completes the proof.
\end{proof}

The amplitude normalization is compatible with the pulled-back
linear potential:
\begin{equation}\label{bg:relative-potential}
\sigma_i^4(W_\vartheta\circ T_i^\vartheta)^{p-1}
          =\overline W_i^{\,p-1},\qquad a(p-1)=4.
\end{equation}
This identity is used when $p-1<1$.

\section{The projected linear equation}\label{an:section}
The parameters fixed above satisfy
\begin{equation}\label{an:parameter-range}
 s\ge1,\qquad 2<\beta<4,\qquad 3-\beta+2\beta/s>0.
\end{equation}
Indeed, $s=1$ gives the last expression $3+\beta>0$; for $s\ge2$
our configurations have $2<\beta<3$. These are the parameter
restrictions used in the analytic estimates below.
We use the caps and weighted space in Section~\ref{bg:section}.
Set $V=Q^{p-1}$, $Z_0=aQ+y\cdot\nabla Q$ and
$Z_j=\partial_jQ$ for $1\le j\le n$. In the Gram constant below,
$A=A_n$. All constants may depend on the fixed geometry and dimension.
The multiplier fields are $Z_{i,\alpha}(X)=Z_\alpha((X-q_i)/\mu_i)$.
The full translation space is convenient for the local inverse; restricting
to the prescribed symmetry leaves only dilation and radial translation.

\subsection{A uniform augmented single-bubble inverse}
For $0<\tau<\min(s,4)$, set
\[
 Y_\tau=\{\psi\in C(\mathbb R^n):
 \|\psi\|_\tau:=\sup_y(1+|y|)^\tau|\psi(y)|<\infty\}.
\]
The choice $\tau=1/2$ works for every $n\ge5$.
Let $\chi_L(y)=\chi(y/L)$, where $\chi$ is a fixed smooth radial cutoff,
$0\le\chi\le1$, equal to one on $B_{1/2}$ and supported in $B_1$.

\begin{theorem}\label{an:single}
There are $L_0,C>0$ such that, for every $L\ge L_0$,
$h\in L^\infty(\mathbb R^n)$ and $d\in\mathbb R^{n+1}$, the system
\begin{align}
 \psi&=\Gamma*\left[\chi_LV
       \left(p(h+\psi)+\sum_{\alpha=0}^n C_\alpha Z_\alpha\right)\right],
       \label{an:single-eq}\\
 \int\chi_LVZ_\alpha(h+\psi)&=d_\alpha\qquad(0\le\alpha\le n)
 \label{an:single-orth}
\end{align}
has a unique solution in $Y_\tau\times\mathbb R^{n+1}$, with
\[
 \|\psi\|_\tau+|C|\le C(\|h\|_\infty+|d|).
\]
The same assertion holds with all cutoffs removed.
\end{theorem}
\begin{proof}
Splitting the integral into $|Y|<|X|/2$, $|X-Y|<|X|/2$, and the remaining
region gives, for $b>4$,
\[
 \int\Gamma(X-Y)(1+|Y|)^{-b}\,dY
 \le C\begin{cases}
 (1+|X|)^{4-b},&4<b<n,\\
 (1+|X|)^{-s}\log(2+|X|),&b=n,\\
 (1+|X|)^{-s},&b>n.
 \end{cases}
\]
Thus $h\mapsto\Gamma*(Vh)$ is bounded from $L^\infty$ to $Y_\tau$.
For $T\psi=\Gamma*(pV\psi)$ and any
$0<\delta<\min(4,s-\tau)$, the same bound yields
$\|T\psi\|_{\tau+\delta}\le C\|\psi\|_\tau$.
On compact sets the sources are uniformly bounded, and the kernel and its
first derivatives are locally integrable. The images are therefore locally
equicontinuous, with uniformly vanishing weighted tails. Arzel\`a--Ascoli
proves compactness of $T$ on $Y_\tau$.

For the kernel classification, put $\Omega=2/(1+|y|^2)$ and write
$\psi=\Omega^av$ under stereographic projection. The conformal identity is
\[
 \Delta^2(\Omega^av)=\Omega^{a+4}P_nv,\qquad
 P_n=\left(-\Delta_{S^n}+\frac{n(n-2)}4\right)
     \left(-\Delta_{S^n}+\frac{(n-4)(n+2)}4\right).
\]
Also $V=\kappa_n\Omega^4$, where $\kappa_n=n(n-4)(n^2-4)/16$.
If $(I-T)\psi=0$, then $(P_n-p\kappa_n)v=0$ off the north pole and
$|v|\le Cr^{-s+\tau}$ at distance $r$ from it. This growth is locally integrable.
Testing with a cutoff at radius $\varepsilon$ and transferring all derivatives
to the test function, the fourth-order cutoff terms are
$O(\varepsilon^{n-4-s+\tau})=O(\varepsilon^\tau)$; lower-order terms are smaller.
The equation extends distributionally to the sphere, and elliptic regularity
makes $v$ smooth. With $\lambda_l=l(l+n-1)$, the eigenvalue of $P_n$ is
\[
 (\lambda_l+n(n-2)/4)(\lambda_l+(n-4)(n+2)/4).
\]
It is strictly increasing in $l$, and equals $p\kappa_n$ precisely for $l=1$.
Consequently $\ker(I-T)=\operatorname{span}\{Z_0,\ldots,Z_n\}$.

For sphere coordinates $\xi_0=(|y|^2-1)/(1+|y|^2)$ and
$\xi_j=2y_j/(1+|y|^2)$, one has
$Z_\alpha=-aA2^{-a}\Omega^a\xi_\alpha$. Hence
\[
 \Gamma*(VZ_\alpha)=Z_\alpha/p,\qquad
 \int VZ_\alpha Z_\gamma=g_n\delta_{\alpha\gamma},\quad
 g_n=\frac{a^2A^{p+1}|S^n|}{2^n(n+1)}>0.
\]
The convolution identity also follows from the Jacobi equation and uniqueness
of an entire biharmonic function tending to zero. The functionals
$\ell_\alpha(f)=\int VZ_\alpha f$ are continuous on $Y_\tau$ and satisfy
$\ell_\alpha(Tf)=\ell_\alpha(f)$. The Fubini interchange is absolute, since
$V|Z_\alpha|=O((1+r)^{-n-4})$.
Fredholm index zero and independence of these functionals, furnished by the
Gram matrix, give
\[
 \operatorname{Ran}(I-T)=\bigcap_{\alpha=0}^n\ker\ell_\alpha.
\]
Applying $\ell_\alpha$ to the untruncated augmented equation gives
$C_\alpha=-p\ell_\alpha(h)/g_n$. Invert on the displayed closed complement,
then add Jacobi components to impose
$\ell_\alpha(\psi)=d_\alpha-\ell_\alpha(h)$.
This proves existence, uniqueness, and the bounded inverse without cutoffs.

For truncation, the same three-region computation gives
\[
 \sup_X(1+|X|)^\tau
 \int_{|Y|>L/2}\Gamma(X-Y)(1+|Y|)^{-8-\tau}\,dY\le CL^{-4}.
\]
For $|X|\le L/4$ the tail is at most
$CL^{-4-\tau}(1+|X|)^\tau$. For the other points, split at the scale $|X|$;
any critical logarithm is absorbed using
$\sup_{r\ge1}r^{\tau-s}\log r<\infty$.
The potential block therefore changes by $O(L^{-4})$ in operator norm.
The multiplier columns obey the same bound because $Z_\alpha\in Y_\tau$.
The orthogonality rows change by at most
\[
 \int_{|Y|>L/2}V|Z_\alpha|(1+|Y|)^{-\tau}\,dY\le CL^{-4-\tau}.
\]
The truncated right-hand side acting on $h$ remains uniformly bounded.
Neumann perturbation of the augmented inverse proves the theorem.
\end{proof}

\subsection{Physical scaling and the two source quantities}
Define
\[
 J_0(L)=1+\int_1^Lr^{s-5}\,dr\qquad(L\ge1).
\]
Scale to $Q_i(X)=\mu_i^{-a}Q((X-q_i)/\mu_i)$ with microscopic cutoff
$L_i\asymp R_i/\mu_i$. Normalize the orthogonality rows by
\[
 \int\chi_iQ_i^{p-1}Z_{i,\alpha}(h_i+\psi_i)
      =\mu_i^sR_i^\beta d_{i,\alpha}.
\]
For $b_i=\|h_i\|_\infty+R_i^\beta|d_i|$, the response source
\[
 f_i=\chi_iQ_i^{p-1}
       \left[p(h_i+\psi_i)+\sum_\alpha C_{i,\alpha}Z_{i,\alpha}\right],
 \qquad\psi_i=\Gamma*f_i
\]
satisfies
\begin{align}
 |\psi_i(X)|&\le Cb_i(1+|X-q_i|/\mu_i)^{-\tau},&|C_i|&\le Cb_i,
 \label{an:local-decay}\\
 \|\Gamma*|f_i|\|_\infty&\le Cb_i,&
 \|f_i\|_1&\le C\mu_i^sJ_0(L_i)b_i.\label{an:local-mass}
\end{align}
Indeed $Q_i^{p-1}\,dX=\mu_i^sV\,dy$, while the potential scaling factor is
$\mu_i^{n-s-4}=1$. The bracket is bounded by $Cb_i$ in microscopic variables.
Boundedness of $\Gamma*V$ and the radial mass integral prove
\eqref{an:local-mass}. The factor $J_0$ cannot be omitted when $n\ge8$.

\subsection{The model inverse for countably many caps}
Use the fixed $t,\beta$ and caps of Section~\ref{bg:section},
and choose $1<\rho<t^{-(4-\beta)}$ as in Proposition~\ref{sym:rightinverse}.
Here $i=(k,\nu)$ denotes a point on the $k$th shell and $R_i=R_k$.
The enlarged caps are disjoint, lie in $D=B_1$, and satisfy $|X|\asymp R_i$.
Write $A_k=\mathcal A_k$ for the fixed annuli in \eqref{bg:space}, retain
the space $\mathcal Y$, and set
\[
 \|C\|_{\mathcal C}=\sum_{k\ge0}\rho^kR_k^{-\beta}\max_\nu|C_{k,\nu}|,
 \quad\|d\|_{\mathcal D}=\sum_{k\ge0}\rho^k\max_\nu|d_{k,\nu}|.
\]
Use the prescribed invariant function space. Multiplier coefficients and
orthogonality data are equivariant: the dilation component transforms as
a scalar and the translation components as vectors. The operators
$N_2$ and $\mathcal S$ are those defined in \eqref{bg:S}.

\begin{theorem}\label{an:countable}
If all $L_i$ are sufficiently large and
\[
 \eta_*:=\sup_i(\mu_i/R_i)^sJ_0(L_i)
\]
is sufficiently small, the model augmented system
\begin{align}
 \phi-\mathcal S\left[
 p\sum_i\chi_iQ_i^{p-1}\phi+
 \sum_{i,\alpha}C_{i,\alpha}\chi_iQ_i^{p-1}Z_{i,\alpha}\right]&=g,
 \label{an:countable-eq}\\
 \frac1{\mu_i^sR_i^\beta}
 \int\chi_iQ_i^{p-1}Z_{i,\alpha}\phi&=d_{i,\alpha}
 \label{an:countable-orth}
\end{align}
has a bounded two-sided inverse, with
$\|\phi\|_{\mathcal Y}+\|C\|_{\mathcal C}
\le C(\|g\|_{\mathcal Y}+\|d\|_{\mathcal D})$.
For $\mu_i\asymp\lambda^{2/s}R_i^{2+2\beta/s}$, these hypotheses hold for
sufficiently small $\lambda$ in every fixed dimension.
\end{theorem}
\begin{proof}
For local inputs $h_i$, solve the local augmented problems from
Theorem~\ref{an:single}. On the Banach space
$\|h\|_{\mathcal H}=\sum_k\rho^kR_k^{-\beta}
\max_\nu\|h_{k,\nu}\|_\infty$, solve
\[
 h_i=g|_{K_i}+\left(\sum_jN_2\psi_j-\psi_i\right)\big|_{K_i}.
\]
An outer source shell $\ell<k$ is first combined as a complete even shell;
its normalized kernel is bounded by $CR_k^4R_\ell^{-s-4}$.
An inner shell contributes $CR_k^2R_\ell^{-s-2}$.
After deleting the single self-response, the target shell is not even as a whole.
Its other caps are macroscopically separated and give $CR_k^{-s}$;
the remaining self term $-T_2\psi_i$ obeys the same mass estimate.
Consequently the normalized interaction matrix is $C\eta_\ell H_{k\ell}$, where
\[
 H_{k\ell}=\begin{cases}
 t^{(4-\beta)(k-\ell)},&k>\ell,\\
 1,&k=\ell,\\
 t^{(\beta-2)(\ell-k)},&k<\ell.
 \end{cases}
\]
The assumptions ensure
$\sup_\ell\rho^{-\ell}\sum_k\rho^kH_{k\ell}<\infty$.
The linear fixed-point map thus has norm at most $C\eta_*<1/2$, and the
orthogonality data contribute at most $C\eta_*\|d\|_{\mathcal D}$.
Its unique fixed point obeys
$\|h\|_{\mathcal H}\le C(\|g\|_{\mathcal Y}+\|d\|_{\mathcal D})$.

Set $\phi=g+\sum_iN_2\psi_i$.
At comparable scales a source contributes at most
$\|\Gamma*|f_i|\|_\infty+R_i^{-s}\|f_i\|_1$;
all other scales use the preceding matrix. This proves the estimates for
$\phi$ and $C$. Writing $b_i=R_i^\beta\xi_i$, one has
\[
 \sum_iR_i^{-s-2}\|f_i\|_1
 \le C\sum_i\eta_iR_i^{\beta-2}\xi_i<\infty.
\]
The absolute second kernel moment converges, so the normalized series equals
the integral and solves \eqref{an:countable-eq}--\eqref{an:countable-orth}.
Conversely, for any global solution define its cap sources $f_i$,
$\psi_i=\Gamma*f_i$, and $h_i=\phi|_{K_i}-\psi_i|_{K_i}$.
The direct source potential bound puts $h$ in $\mathcal H$.
It solves the same local systems and fixed-point equation. Both uniqueness
statements therefore imply global uniqueness, proving a two-sided inverse.

Finally,
\[
 \eta_i\le C\begin{cases}
 (\mu_i/R_i)^s,&s<4,\\
 (\mu_i/R_i)^4[1+\log(R_i/\mu_i)],&s=4,\\
 (\mu_i/R_i)^4,&s>4.
 \end{cases}
\]
The matching scales give $\eta_*\le C\lambda^2$ for $s=1,2,3$,
$C\lambda^2(1+|\log\lambda|)$ for $s=4$, and $C\lambda^{8/s}$ for $s>4$.
Also $\inf_iL_i\to\infty$. This proves the last assertion.
\end{proof}

\subsection{Actual parameter differences in moving coordinates}
Let $\varepsilon=\|\vartheta-\vartheta'\|_X$, and use the maps
$\Phi_\vartheta$ and $T_i^\vartheta$ from \eqref{bg:T}--\eqref{bg:Phi}.
Write $\mu_i^0=\mu_i(0)$, $q_i^0=q_i(0)$,
$K_i^0=B(q_i^0,\delta R_i)$ and
$\sigma_i^\vartheta=\mu_i^\vartheta/\mu_i^0$.
These maps preserve the annuli and agree with $T_i^\vartheta$ on
$B(q_i^0,10\delta R_i)$. In addition to \eqref{bg:Phi}, they satisfy
$|D^2\Phi_\vartheta(X)|\le C/|X|$. Shrink the parameter ball so that
$\|D\Phi_\vartheta-I\|_\infty<1/2$; straight interpolations between two
maps then have uniform bi-Lipschitz constants.

Denote the model operator in
\eqref{an:countable-eq}--\eqref{an:countable-orth} by
$\mathcal A_\vartheta^{\rm mod}$, and set
\begin{equation}\label{an:pull-operator}
 \Pi_\vartheta v=v\circ\Phi_\vartheta,\qquad
 \widetilde{\mathcal A}_\vartheta
 =(\Pi_\vartheta\oplus I_{\mathcal D})
   \mathcal A_\vartheta^{\rm mod}
   (\Pi_\vartheta^{-1}\oplus I_{\mathcal C}).
\end{equation}
Thus both $\phi$ and $g$ are pulled back by composition alone.
The coefficients $C$ and the normalized orthogonality data $d$ are
unchanged. In particular, no amplitude factor is included in
$\widetilde\phi=\Pi_\vartheta\phi$.

Put $w_i^0=\chi_i^0(Q_i^0)^{p-1}$ and
$Z_{i,\alpha}^0(Y)=Z_\alpha((Y-q_i^0)/\mu_i^0)$.
Changing variables on the source cap gives the effective kernel
\begin{equation}\label{an:pull-kernel}
 \mathscr K_i^\vartheta(X,Y)
   =(\sigma_i^\vartheta)^s
          K_2(\Phi_\vartheta X,T_i^\vartheta Y),
 \qquad Y\in K_i^0.
\end{equation}
Indeed the potential and volume element contribute
$(\sigma_i^\vartheta)^{-4}$ and $(\sigma_i^\vartheta)^n$, respectively.
Consequently
\begin{align}
 \widetilde{\mathcal A}_\vartheta(\widetilde\phi,C)
   &=\bigl(\widetilde\phi-\widetilde{\mathcal T}_\vartheta
                (\widetilde\phi,C),\,\ell^0(\widetilde\phi)\bigr),
       \label{an:pull-block}\\
 \widetilde{\mathcal T}_\vartheta(\widetilde\phi,C)(X)
   &=\sum_i\int_{K_i^0}\mathscr K_i^\vartheta(X,Y)F_i(Y)\,dY,
       \notag\\
 F_i&=w_i^0\left(p\widetilde\phi+
                 \sum_\alpha C_{i,\alpha}Z_{i,\alpha}^0\right),
 \qquad
 \ell_{i,\alpha}^0(\widetilde\phi)
   =\frac{\int_{K_i^0}w_i^0Z_{i,\alpha}^0\widetilde\phi}
          {(\mu_i^0)^sR_i^\beta}.
       \notag
\end{align}
The last row is independent of the parameters because
$(\mu_i^\vartheta)^{-s}(\sigma_i^\vartheta)^s=(\mu_i^0)^{-s}$.

\begin{proposition}\label{an:param}
In the coordinates \eqref{an:pull-operator},
\[
 \|\widetilde{\mathcal A}_{\vartheta}
       -\widetilde{\mathcal A}_{\vartheta'}\|
 +\|\widetilde{\mathcal A}_{\vartheta}^{-1}
       -\widetilde{\mathcal A}_{\vartheta'}^{-1}\|
 \le C\varepsilon\eta_*.
\]
The norms are the respective operator norms between
$\mathcal Y\times\mathcal C$ and $\mathcal Y\times\mathcal D$,
and $\eta_*$ is the supremum in Theorem~\ref{an:countable}
over the entire parameter ball.
\end{proposition}
\begin{proof}
For an arbitrary input $(\widetilde\phi,C)$, define
\[
 b_i^{\rm in}=\|\widetilde\phi\|_{L^\infty(K_i^0)}+|C_i|,
 \qquad
 \eta_i^0=(\mu_i^0/R_i)^sJ_0(L_i).
\]
The Jacobi fields are bounded. Direct integration of $w_i^0$ therefore
gives, without solving a local equation,
\begin{equation}\label{an:pull-input-mass}
 \|F_i\|_1\le C(\mu_i^0)^sJ_0(L_i)b_i^{\rm in}
             =C\eta_i^0R_i^s b_i^{\rm in}.
\end{equation}

For $X\in B(q_i^0,10\delta R_i)$ and $Y\in K_i^0$, the singular part
of \eqref{an:pull-kernel} is exactly
\begin{equation}\label{an:pull-exact}
 (\sigma_i^\vartheta)^s
     \Gamma(T_i^\vartheta X-T_i^\vartheta Y)=\Gamma(X-Y).
\end{equation}
This fixes the same-cap uncompensated principal block and its multiplier
columns. For the Taylor subtraction, put
\[
 \mathscr K_{i,\mathrm{sub}}^\vartheta(X,Y)
 =(\sigma_i^\vartheta)^s
   \left[T_{2,A}\Gamma(A-T_i^\vartheta Y)\right]_{A=T_i^\vartheta X},
\]
where $T_{2,A}$ is taken at $A=0$. The coefficients of degree $r\le2$
are $O(R_i^{-s-r})$, and their parameter differences have the same
bound multiplied by $\varepsilon$. Since
$|\delta q_i|\le C\varepsilon R_i$ and
$|\delta\sigma_i|\le C\varepsilon$, differentiation in $X$ yields
\begin{equation}\label{an:pull-self}
 |D_X^j\delta\mathscr K_{i,\mathrm{sub}}(X,Y)|
       \le C\varepsilon R_i^{-s-j},\qquad 0\le j\le2.
\end{equation}
Here and below $\delta$ denotes the difference between $\vartheta$
and $\vartheta'$. Equations \eqref{an:pull-input-mass}--\eqref{an:pull-self}
give
\begin{equation}\label{an:pull-self-source}
 \left|D_X^j\int_{K_i^0}\delta\mathscr K_i(X,Y)F_i(Y)\,dY\right|
       \le C\varepsilon\eta_i^0R_i^{-j}b_i^{\rm in}.
\end{equation}
Thus the full multiplier columns have small parameter differences;
only their same-cap singular principal parts are exactly fixed.

For completeness, the complete even kernel obeys the following bounds.
Set
$\mathscr G_\vartheta(X,Y)
=K_{2,e}(\Phi_\vartheta X,\Phi_\vartheta Y)\det D\Phi_\vartheta(Y)$.
For $0\le j\le2$,
\begin{equation}\label{an:pull-three-regions}
 |D_X^j\delta\mathscr G(X,Y)|
 \le C\varepsilon
 \begin{cases}
 |X|^{4-j}|Y|^{-s-4},&|X|\le |Y|/4,\\
 |X|^{2-j}|Y|^{-s-2},&|Y|\le |X|/4,\\
 |X|^{-s-j}+|X-Y|^{-s-j}+|X+Y|^{-s-j},
       &|X|/4<|Y|<4|X|.
 \end{cases}
\end{equation}
To verify the first two cases, use the even Taylor remainder and
$D_AK_{2,e}=O(|A|^3|B|^{-s-4})$ for outer sources,
or $O(|A||B|^{-s-2})$ for inner sources.
The chain term $D_AK_{2,e}D^2\Phi_\vartheta$ and its difference
then have the stated second-derivative order.
At comparable scales, oddness and the uniform bi-Lipschitz bounds give
$|\delta\Phi(X)\mp\delta\Phi(Y)|\le C\varepsilon|X\mp Y|$.
The fundamental-kernel derivatives, the Taylor polynomial and the
Jacobian can now be estimated separately. This proves
\eqref{an:pull-three-regions}. The same far-region bounds apply to
\eqref{an:pull-kernel} after pairing opposite source caps; on a source
cap its additional factor relative to $\mathscr G_\vartheta$ is
$(\sigma_i^\vartheta)^{-4}$, whose size is bounded and whose difference
is $O(\varepsilon)$.

Choose a fixed integer $J$ so that $|k-\ell|>J$ puts every observation
point in $A_k$ and every source cap on shell $\ell$ in a far region.
For $|k-\ell|\le J$, a source and observation point are either in the
same affine cap, where \eqref{an:pull-self-source} applies, or separated
by at least $cR_\ell$. Indeed the source lies in the $\delta R_\ell$
ball and the affine patch contains the $10\delta R_\ell$ ball;
distinct caps satisfy the fixed separation condition.
The latter case has kernel difference at most
$C\varepsilon R_\ell^{-s}$.
There are at most $N(2J+1)$ such source caps for each target shell,
and the finitely many radius ratios are absorbed in the constant.
Complete distant shells are even by equivariance.
Their contributions, using \eqref{an:pull-input-mass} and
\eqref{an:pull-three-regions}, consequently satisfy
\begin{equation}\label{an:pull-shell-bound}
 R_k^{-\beta}
 \|\delta\widetilde{\mathcal T}(\widetilde\phi,C)\|_{L^\infty(A_k)}
 \le C\varepsilon\sum_{\ell\ge0}
       \eta_\ell^0H_{k\ell}\,\xi_\ell,
 \qquad
 \xi_\ell=R_\ell^{-\beta}\max_{i\text{ on shell }\ell}b_i^{\rm in},
\end{equation}
where $\eta_\ell^0=\max_{i\text{ on shell }\ell}\eta_i^0$ and
$H_{k\ell}$ is the matrix in Theorem~\ref{an:countable}.
The weighted column sum is explicitly bounded by
\begin{equation}\label{an:pull-column-sum}
 \sup_{\ell\ge0}\sum_{k\ge0}\rho^{k-\ell}H_{k\ell}
 \le 1+
 \frac{\rho t^{4-\beta}}{1-\rho t^{4-\beta}}
 +\frac{t^{\beta-2}/\rho}{1-t^{\beta-2}/\rho}<\infty.
\end{equation}
Moreover
$\sum_\ell\rho^\ell\xi_\ell
\le C(\|\widetilde\phi\|_{\mathcal Y}+\|C\|_{\mathcal C})$.
Summing \eqref{an:pull-shell-bound} therefore proves the operator
difference bound. The second row has zero difference by
\eqref{an:pull-block}. Theorem~\ref{an:countable} and the uniform
boundedness of the pullbacks give uniform two-sided inverse bounds.
Finally,
\[
 \widetilde{\mathcal A}_\vartheta^{-1}
 -\widetilde{\mathcal A}_{\vartheta'}^{-1}
 =\widetilde{\mathcal A}_\vartheta^{-1}
   (\widetilde{\mathcal A}_{\vartheta'}
           -\widetilde{\mathcal A}_\vartheta)
   \widetilde{\mathcal A}_{\vartheta'}^{-1},
\]
which proves the asserted inverse estimate.
\end{proof}

\section{The nonlinear projected solution}\label{anq:section}
We retain the notation of the preceding section. In a cap, write
$R=R_i$, $d=|X-q_i|$ and $L=R/\mu_i$. By Lemma~\ref{bg:bounds},
$\mu_i^a\asymp\lambda R_i^{s+\beta}$.
The weight $\rho>1$ is chosen sufficiently close to one, on a zero-free
circle from Proposition~\ref{sym:rightinverse}, that every positive shell
power below is summable. The number of such restrictions is finite.

\subsection{Separate control of potential amplitude and source mass}
\begin{lemma}\label{anq:two-quantities}
Suppose the sources $f_i$, supported in their respective caps, form an
equivariant family under $\mathfrak G$. Set
\[
 A_i=\|\Gamma*|f_i|\|_\infty,\qquad M_i=\|f_i\|_1.
\]
Then
\[
 \left\|\mathcal S\sum_i f_i\right\|_{\mathcal Y}
 \le C\sum_{k\ge0}\rho^kR_k^{-\beta}
             \max_{i\text{ on shell }k}(A_i+R_i^{-s}M_i).
\]
For a complete even annular source satisfying $|f|\le F_k$ on $A_k$, the
right-hand side can be replaced by $C\sum_k\rho^kR_k^{4-\beta}F_k$.
\end{lemma}
\begin{proof}
For a complete even outer shell, the compensated kernel is bounded by
$Cr^4R^{-s-4}$; for an inner shell, the bound is $Cr^2R^{-s-2}$.
At comparable scales use $A_i+R_i^{-s}M_i$. After normalization, the
interaction matrix is
\[
 H_{k\ell}=\begin{cases}
 t^{(4-\beta)(k-\ell)},&k>\ell,\\
 1,&k=\ell,\\
 t^{(\beta-2)(\ell-k)},&k<\ell.
 \end{cases}
\]
Its weighted column sums are uniformly bounded. For an annular source, the
mass and local potential are bounded by $CR^{s+4}F_k$ and $CR^4F_k$,
respectively. This proves the second assertion.
\end{proof}

\begin{lemma}\label{anq:moments}
Set
\[
 J_j(L)=1+\int_1^L r^{s+j-5}\,dr.
\]
For $j\ge0$,
\[
 \int_{d<R}d^jQ_\mu^{p-1}\le C\mu^{s+j}J_j(L).
\]
The function $J_j$ is comparable to $1$, $1+\log L$, or $L^{s+j-4}$,
according as $s+j<4$, $s+j=4$, or $s+j>4$.
For $0\le j\le2$, one also has
\[
 \|\Gamma*(\mathbf1_{d<R}Q_\mu^{p-1}d^j)\|_\infty\le C\mu^j.
\]
Consequently, if $|B|\le b_0+b_1d+\lambda d^2$, the two quantities of the
source $\mathbf1_{d<R}Q_\mu^{p-1}B$ satisfy
\begin{align*}
 A&\le C(b_0+\mu b_1+\lambda\mu^2),\\
 M&\le C[\mu^sJ_0b_0+\mu^{s+1}J_1b_1+
                      \lambda\mu^{s+2}J_2].
\end{align*}
\end{lemma}
\begin{proof}
Since $Q_\mu^{p-1}=C\mu^4(\mu^2+d^2)^{-4}$, radial scaling gives the mass
formula. For the potential amplitude, scaling reduces the assertion to
boundedness of $I_4[(1+|y|^2)^{-4}|y|^j]$. Split the integral into inner,
comparable and outer regions. The far radial integral is
$\int r^{j-5}\,dr$, which converges for $j<4$, and the local kernel is
integrable. These estimates are uniform in the observation point.
\end{proof}

\begin{remark}
The two quantities should be kept separate. Dividing the mass by $\mu^s$
and then treating it as the potential amplitude introduces unnecessary
logarithms or negative shell powers when $n\ge8$.
Relative to the $j$th potential-amplitude bound, the extra factor in $M/R^s$
is $L^{-s}J_j(L)$. For $j\le2$, it is bounded and tends to zero, including
the critical factor $L^{-2}\log L$ when $s=j=2$.
\end{remark}

\subsection{Only two source derivatives are needed for the bubble tails}
\begin{lemma}\label{anq:tails}
Let $\chi_i$ be macroscopic cap cutoffs with uniform scaling bounds.
Suppose that the complete tail source $h=\sum_i(1-\chi_i)Q_i^p$ is even.
Then $h$ extends through the origin as a $C^2$ function, with norm
$O(\lambda^p)$. In a fixed neighbourhood of the origin,
$V=\Gamma*(\mathbf1_Dh)$ satisfies
\[
 |D^j(V-T_2V)(X)|\le C\lambda^p|X|^{4-j},\qquad j=0,1,2.
\]
If the physical-coordinate parameter differences of the tail sources satisfy
the same termwise derivative bounds multiplied by $\varepsilon$, then the
parameter differences of these estimates have the same extra factor.
\end{lemma}
\begin{proof}
For $j\le2$, termwise scaling gives
\[
 \|D^j[(1-\chi_i)Q_i^p]\|_\infty
       \le C\lambda^pR_i^{\beta p-j}.
\]
Since $\beta p>2$, the series and its first two derivatives converge uniformly.
The local potential has $C^0$ norm $O(\lambda^p)$ and satisfies
$\Delta^2V=h$. Interior Schauder estimates yield
$\|V\|_{C^{4,\alpha}}\le C\lambda^p$ for any fixed $0<\alpha<1$.
Contributions separated from the origin are smooth.
Evenness of $V$ makes its third derivatives at the origin vanish, and Taylor's
formula gives the claimed bounds. Apply the same argument to actual parameter
differences for the final assertion.
\end{proof}
This lemma does not require $\beta p>4$, nor does it assert that the tail
source itself extends as a $C^4$ function in high dimensions.

\subsection{Real-power nonlinearities on a ball preserving positivity}
Lemma~\ref{bg:bounds} gives
\begin{equation}\label{anq:bg-hypotheses}
 W(X)\ge c\lambda|X|^\beta,\qquad
 W\le C(Q_i+\lambda R_i^2)\ \text{on caps},\qquad
 W\le C\lambda|X|^2\ \text{off caps}.
\end{equation}
\begin{proposition}[The background residual]\label{anq:residual}
Let $B_i=W-Q_i$ satisfy the second-order envelope in
Lemma~\ref{anq:moments}, with
\begin{align*}
 b_{0,i}&\le C\lambda R_i^\beta e_i+C\lambda^{1+4/s}R_i^4,\\
 b_{1,i}&\le C\lambda R_i^{\beta-1}e_i+C\lambda^{1+4/s}R_i^3,
 \qquad e_i\ge0.
\end{align*}
If $\|e\|_{\rho}:=\sum_k\rho^k\max_{i\text{ on shell }k}e_i\le1$
and the complete bubble tails satisfy Lemma~\ref{anq:tails}, then
\[
 \left\|\mathcal S\left(W^p-\sum_iQ_i^p\right)\right\|_{\mathcal Y}
 \le C\lambda(\|e\|_\rho+\lambda^{4/s}).
\]
\end{proposition}
\begin{proof}
For $W=Q_i+B_i>0$, the real-power mean value inequality gives
$|W^p-Q_i^p|\le C(Q_i^{p-1}|B_i|+|B_i|^p)$.
Lemmas~\ref{anq:moments} and \ref{anq:two-quantities} control the first term.
After division by $\lambda R_i^\beta$, the noncentral coefficients contain
$\lambda^{4/s}R_i^{4-\beta}$, $\mu_i/R_i$, and
$\mu_i^2R_i^{-\beta}=O(\lambda^{4/s}R_i^{4-\beta+4\beta/s})$,
so all the required weighted sums converge. The second term and the source
off the caps are bounded by the annular-source estimate for
$C\lambda^p r^{2p}$. In the identity
$W^p-\sum Q_i^p=\sum\chi_i(W^p-Q_i^p)+(1-\sum\chi_i)W^p
-\sum(1-\chi_i)Q_i^p$, the final term is controlled by the tail lemma.
Both pure-tail contributions are $O(\lambda^p)$, and $p=1+8/s>1+4/s$.
\end{proof}

Put
\[
 \mathscr N_W(\phi)=(W+\phi)^p-W^p-pW^{p-1}\phi,
\qquad
 K_s(\lambda)=\begin{cases}
 \lambda,&s=1,2,3,\\
 \lambda(1+|\log\lambda|),&s=4,\\
 \lambda^{p-2},&s>4.
 \end{cases}
\]
\begin{proposition}\label{anq:nonlinear}
Let $\|\phi\|_{\mathcal Y},\|\psi\|_{\mathcal Y}\le d_0\le c_0\lambda$,
where $c_0$ is small enough that $|\phi|,|\psi|\le W/2$. Then
\begin{align*}
 \|\mathcal S\mathscr N_W(\phi)\|_{\mathcal Y}
    &\le CK_s(\lambda)d_0^2,\\
 \|\mathcal S[\mathscr N_W(\phi)-\mathscr N_W(\psi)]\|_{\mathcal Y}
    &\le CK_s(\lambda)d_0\|\phi-\psi\|_{\mathcal Y}.
\end{align*}
These assertions hold in every fixed dimension $n\ge5$, including $1<p<2$.
\end{proposition}
\begin{proof}
Taylor's integral formula on the positive interval bounds the remainder by
$CW^{p-2}\phi^2$ and its difference by
$CW^{p-2}(|\phi|+|\psi|)|\phi-\psi|$.
First suppose $p\ge2$, equivalently $s\le8$. Direct radial integration gives
\[
 \|\Gamma*(\mathbf1_{d<R}Q_\mu^{p-2})\|_\infty
 +R^{-s}\int_{d<R}Q_\mu^{p-2}
 \le C\mu^a H_s(L),
\]
where $H_s(L)$ is $1$, $1+\log L$, or $L^{s-4}$, according as
$s<4$, $s=4$, or $4<s\le8$.
The potential integral scales as
$\mu^a\int_0^Lr^3(1+r^2)^{-(4-a)}dr$; the analogous radial mass integral,
after division by $R^s$, is no larger than the displayed bound.
Multiplication by $R^\beta$ gives the respective normalized coefficients
\[
 \lambda R^{s+2\beta},\qquad
 \lambda R^{4+2\beta}(1+\log L),\qquad
 \lambda^{p-2}R^{4+\beta(p-1)}.
\]
The background term uses the annular-source bound and has shell exponent
$2p+\beta>0$. Positive shell powers absorb $\log R$, and
$\log L\le C(1+|\log\lambda|+|\log R|)$, giving the stated coefficient.

When $1<p<2$, use the lower bound for the background instead:
$W^{p-2}\phi^2\le C\lambda^{p-2}r^{\beta(p-2)}\phi^2$.
The normalized exponent in the annular-source estimate is
$4+\beta(p-1)>0$, which gives the same assertion.
Applying these estimates to the remainder difference proves the second bound.
\end{proof}

\begin{proposition}[Actual differences of the background parameter]
\label{anq:nonlinear-parameter}
Suppose $W'$ also satisfies \eqref{anq:bg-hypotheses}, and in common
coordinates $W'\asymp W$ and $|W'-W|\le C\varepsilon W$.
On a common positivity ball,
\[
 \|\mathcal S[\mathscr N_W(\phi)-\mathscr N_{W'}(\psi)]\|_{\mathcal Y}
 \le CK_s(\lambda)
       [d_0\|\phi-\psi\|_{\mathcal Y}+\varepsilon d_0^2].
\]
The same assertion holds when the kernel also changes, provided its
parameter difference satisfies the three-region kernel bounds of this
section multiplied by $\varepsilon$.
\end{proposition}
\begin{proof}
Differentiate Taylor's integral formula in the positive background variable.
For every $p>1$, the derivative is bounded by $CW^{p-3}\phi^2$, so background
variation contributes at most $C\varepsilon W^{p-2}\phi^2$.
The proof of Proposition~\ref{anq:nonlinear} controls this term, while its
second estimate handles variation of $\phi$. The same two-quantity estimate
handles kernel differences.
\end{proof}

\subsection{The actual linear potential}
Write $q=p-1$ and $E_s(\lambda)=\lambda K_s(\lambda)$.
Estimate \eqref{anq:bg-hypotheses}, together with
$|B_i|=|W-Q_i|\le C\lambda R_i^2$, gives
\begin{equation}\label{anq:potential-absolute}
 \left\|\mathcal S\left[
   \left(W^{p-1}-\sum_i\chi_iQ_i^{p-1}\right)\phi\right]\right\|_{\mathcal Y}
 \le CE_s(\lambda)\|\phi\|_{\mathcal Y}.
\end{equation}
The truncated bubble potential in the transition region is an ordinary
annular source and obeys the same bound.
For $q>1$, use
$|W^q-Q_i^q|\le C(\lambda R_i^2Q_i^{q-1}+(\lambda R_i^2)^q)$
and the preceding two-quantity integrals.
For $q\le1$, use $|W^q-Q_i^q|\le C|B_i|^q$ and the annular-source estimate.

For $q<1$, fix the principal bubble in common microscopic coordinates.
Lemma~\ref{bg:bounds} gives
\begin{equation}\label{anq:fine-background}
 |B_i-B_i'|\le C\varepsilon\lambda(R_i^\beta+|z|^2),\qquad
 W,W'\ge c(Q_i+\lambda |z|^2),\qquad Q_i\ge c\lambda R_i^\beta.
\end{equation}
Here $z$ is the normal coordinate to the source subspace; affine scaling factors
have already been included in the coefficients.
The mean value formula and monotonicity of negative powers yield
\begin{align*}
 |(W^q-Q_i^q)-((W')^q-Q_i^q)|
 &\le C\varepsilon[\lambda R_i^\beta Q_i^{q-1}+(\lambda |z|^2)^q]\\
 &\le C\varepsilon\lambda^q(R_i^{\beta q}+|z|^{2q})
 \le C\varepsilon\lambda^qR_i^{2q}.
\end{align*}
The second term is zero at $z=0$. This proves the operator parameter-difference
version of \eqref{anq:potential-absolute}, with an extra factor $\varepsilon$.
For $q>1$, the mean value formula and the coarse background difference give
the same result directly; $q=1$ is linear.
A $q$-H\"older inequality alone would give $\varepsilon^q$, which cannot
replace this argument.

Condition \eqref{anq:fine-background} also gives the relative difference
$|W-W'|\le C\varepsilon W$ inside caps. Off caps one needs the corresponding
lower bound $W\ge c\lambda(r^\beta+|z|^2)$ and a difference of the same order.
All these bounds follow from Lemma~\ref{bg:bounds}, after multiplying
the pulled-back background on cap $i$ by $\sigma_i^a$ to fix $Q_i$.
The identity $a(p-1)=4$ makes this normalization agree with the
$\sigma_i^4$ factor in the pulled-back potential.

\subsection{Existence and parameter dependence}

\begin{proposition}\label{an:projection}
There are $\varepsilon_0,\lambda_0>0$ such that, for
$\|\vartheta\|_X\le\varepsilon_0$ and $0<\lambda<\lambda_0$, there is a unique
small pair $(\phi_\vartheta,C_\vartheta)$ satisfying
\begin{align}
 \phi_\vartheta&=\mathcal S\left[(W_\vartheta+\phi_\vartheta)^p
          -\sum_iQ_i^p+\mathcal E_\vartheta\right],\label{an:projected-integral}\\
 \int\chi_i Q_i^{p-1}Z_{i,\alpha}\phi_\vartheta&=0,
 &\mathcal E_\vartheta&=\sum_{i,\alpha}C_{i,\alpha}\chi_iQ_i^{p-1}Z_{i,\alpha}.
 \label{an:projected-orth}
\end{align}
Here the sum runs over all inner centres. Moreover, $U_\vartheta=W_\vartheta+\phi_\vartheta>0$,
and
\begin{align}
 \|\phi_\vartheta\|_{\mathcal Y}+\|C_\vartheta\|_{\mathcal C}
 &\le C\lambda(\|\vartheta\|_X+\lambda^{4/s}),\label{an:projected-bound}\\
 \|\widetilde\phi_\vartheta-\widetilde\phi_{\vartheta'}\|_{\mathcal Y}
       +\|C_\vartheta-C_{\vartheta'}\|_{\mathcal C}
 &\le C\lambda\|\vartheta-\vartheta'\|_X.\label{an:projected-lip}
\end{align}
Tildes denote pullback by the maps in Section~\ref{bg:section}.
The function $U_\vartheta$ is smooth away from the origin and satisfies
$\Delta^2U_\vartheta=U_\vartheta^p+\mathcal E_\vartheta$ there.
\end{proposition}
\begin{proof}
Put $E_s(\lambda)=\lambda K_s(\lambda)$, and denote the actual augmented
operator by
\[
 \mathcal A_\vartheta(\phi,C)=
 \left(\phi-\mathcal S[pW_\vartheta^{p-1}\phi+\mathcal E_C],
 \left(\frac{\int\chi_iQ_i^{p-1}Z_{i,\alpha}\phi}
 {\mu_i^sR_i^\beta}\right)_{i,\alpha}\right).
\]
By \eqref{anq:potential-absolute}, this is an $O(E_s(\lambda))$ perturbation
of the model in Theorem~\ref{an:countable}. It has a uniformly bounded
inverse. Proposition~\ref{an:param} and the actual-potential difference
bound give an inverse difference of size
$C(\eta_*+E_s(\lambda))\|\vartheta-\vartheta'\|_X$ in common coordinates.

Let $r_\vartheta=\mathcal S(W_\vartheta^p-\sum_iQ_i^p)$. Lemma~\ref{bg:bounds},
Proposition~\ref{sym:regularity}, and Proposition~\ref{anq:residual} imply
\begin{equation}\label{an:residual-lip}
 \|r_\vartheta\|_{\mathcal Y}\le C\lambda(\|\vartheta\|_X+\lambda^{4/s}),
 \qquad
 \|\widetilde r_\vartheta-\widetilde r_{\vartheta'}\|_{\mathcal Y}
       \le C\lambda\|\vartheta-\vartheta'\|_X.
\end{equation}
For completeness, the difference estimate uses the actual differences of
the central value and gradient of $B_i$, controlled in weighted $\ell^1$
by the point-source map. The quadratic remainder contributes
$C\lambda\varepsilon\mu_i^2/R_i^\beta$, where
$\varepsilon=\|\vartheta-\vartheta'\|_X$. This is summable because
\[
 \frac{\mu_i^2}{R_i^\beta}
 \le C\lambda^{4/s}R_i^{4-\beta+4\beta/s}.
\]
The pure background difference is bounded by
$C\varepsilon\lambda^p r^{2p}$. The complete tail source and its actual
parameter difference have $C^2$ norm at most $C\lambda^p(1+\varepsilon)$,
by Lemma~\ref{anq:tails}. The two-quantity potential bound and the kernel
difference bound in Proposition~\ref{an:param} control these terms.
Thus \eqref{an:residual-lip} does not follow from formal differentiation
of an absolute residual estimate.

Solve
\[
 (\phi,C)=\mathcal A_\vartheta^{-1}
       (r_\vartheta+\mathcal S\mathscr N_{W_\vartheta}(\phi),0).
\]
On a ball of radius $C\lambda(\varepsilon_0+\lambda^{4/s})$, the contraction
constant is at most $CE_s(\lambda)(\varepsilon_0+\lambda^{4/s})$.
First choose $\varepsilon_0$ small enough to preserve $W+\phi>0$ using
$W\ge c\lambda r^\beta$, and then choose $\lambda_0$ small. The contraction
principle gives the solution and \eqref{an:projected-bound}.
Subtracting the two fixed-point equations, using
\eqref{an:residual-lip}, the inverse difference, and
Proposition~\ref{anq:nonlinear-parameter}, and absorbing the contraction
term proves \eqref{an:projected-lip}. This argument does not use the
effective background or a Pohozaev identity.
The source has a finite absolute second kernel moment by the preceding
estimates. Hence the integral equation implies the asserted equation
distributionally. On compact subsets of the punctured ball, its right-hand
side is locally bounded; elliptic regularity, positivity, and iteration
give smoothness.
\end{proof}

\section{The effective background and the reduced equations}\label{an:effective-section}

Write $i=(k,\nu)$ for an inner centre and let $\zeta_i$ be the transported
radial cutoff equal to one on $\Omega_i(2)$ and supported in $\Omega_i(3)$.
The multiplier cutoff $\chi_i$ is supported in $\Omega_i(1)$.
Define
\begin{align}
 \mathcal E_i&=\sum_\alpha C_{i,\alpha}\chi_iQ_i^{p-1}Z_{i,\alpha},&
 f_i&=\zeta_iU_\vartheta^p+\mathcal E_i,\notag\\
 V_i&=\Gamma*f_i,& H_i&=U_\vartheta-V_i,\label{an:source-definition}\\
 \delta f_i&=\zeta_i(U_\vartheta^p-Q_i^p)+\mathcal E_i,&
 f_{\rm out}&=(1-\sum_i\zeta_i)U_\vartheta^p-\sum_i(1-\zeta_i)Q_i^p.
 \label{an:delta-source}
\end{align}
Except for an explicitly stated whole-space bubble tail, sources are extended
by zero outside $B_1$. The function $H_i$ is biharmonic on $\Omega_i(2)$,
and $\Delta^2H_i=(1-\zeta_i)U_\vartheta^p$ on $\Omega_i(4)$.

\subsection{Masses and source differences}

\begin{lemma}\label{an:moments}
Let $m_i=\int f_i$ and $I_p=\int_{\R^n}Q^p$. Uniformly in the inner centres,
\begin{equation}\label{an:moment-estimate}
 m_i=I_p\mu_i^a[1+o(1)]>0,
 \qquad \int |X-q_i|^j|f_i(X)|\dd X\le C\mu_i^{a+j},
 \quad j=0,1,2,
\end{equation}
as $\lambda\downarrow0$. For
$\widehat f_i(y)=\mu_i^{a+4}f_i(q_i+\mu_i y)$,
\begin{equation}\label{an:micro-source-lip}
 \|\widehat f_i^\vartheta-\widehat f_i^{\vartheta'}\|_{L^1((1+|y|)\dd y)}
       \le\omega(\lambda)\|\vartheta-\vartheta'\|_X,
 \qquad \omega(\lambda)\longrightarrow0.
\end{equation}
\end{lemma}
\begin{proof}
On the source support, $w_i=U_\vartheta-Q_i$ satisfies
$|w_i|\le C\lambda R_i^2$, and its common-coordinate difference satisfies
the same bound multiplied by $\varepsilon=\|\vartheta-\vartheta'\|_X$.
Also $|C_i|\le C\lambda R_i^\beta$ with the same difference bound.
Use
\[
 |(Q_i+w_i)^p-Q_i^p|
 \le C(Q_i^{p-1}|w_i|+|w_i|^p).
\]
The $j$th moment of $Q_i^p$ is $O(\mu_i^{a+j})$, and the lost mass outside
the cutoff is relatively $O((\mu_i/R_i)^4)$. With $L=R_i/\mu_i$, the
relative error from the first term is bounded by
$C\lambda R_i^2\mu_i^aJ_j(L)$. For $j=2$ this is at most
\[
 C\begin{cases}
 \lambda^2R_i^{\beta+3},&s=1,\\
 \lambda^2R_i^{\beta+4}(1+\log L),&s=2,\\
 \lambda^{4/s}R_i^{4-\beta+4\beta/s},&s>2.
 \end{cases}
\]
The corresponding $j=0,1$ errors also tend uniformly to zero.
The pure background relative $j$th moment is bounded by
\[
 C\lambda^{(8-2j)/s}
 R_i^{4+2p-\beta-j-2\beta j/s}.
\]
For $j=2$, the shell exponent is $(4-\beta)(1+4/s)>0$;
the other two exponents are larger. The multiplier moments are bounded
by $C\lambda R_i^\beta\mu_i^a$, since
$\int |y|^jV|Z_\alpha|\dd y<\infty$ for $j\le2$.
These estimates prove \eqref{an:moment-estimate}.

The microscopic cutoff is independent of the parameter, and exactly
\[
 \widehat f_i=\zeta_{L_i}[Q+\mu_i^aw_i(q_i+\mu_i y)]^p
       +\mu_i^a\sum_\alpha C_{i,\alpha}\chi_{L_i}VZ_\alpha.
\]
The positive-interval mean value inequality bounds its weighted actual
difference by $C\varepsilon$ times
\begin{equation}\label{an:source-lip-terms}
 \lambda\mu_i^aR_i^2[J_0(L_i)+J_1(L_i)]
 + (\lambda\mu_i^aR_i^2)^pL_i^{s+5}
 +\lambda^2R_i^{s+2\beta}.
\end{equation}
The second term equals, up to a fixed factor,
$\lambda^{6/s}R_i^{5-\beta+(16-2\beta)/s}$.
For the first term, $s=1$ gives $O(\lambda^2R_i^{\beta+3})$,
$s=2$ gives $O(\lambda^2R_i^{\beta+4})$,
$s=3$ gives $O(\lambda^2R_i^{\beta+5}(1+\log L_i))$, and $s>3$ gives
$O(\lambda^{6/s}R_i^{5-\beta+6\beta/s})$.
All exponents of $R_i$ are positive, so the expression tends uniformly
to zero. This proves \eqref{an:micro-source-lip}, including the difference
of the normalizing mass $m_i/\mu_i^a=\int\widehat f_i$.
\end{proof}

\subsection{The effective background estimates}

\begin{proposition}\label{an:background}
The sequences
\[
 \left(\frac{H_i(q_i)}{\lambda R_i^\beta},
       \frac{q_i\cdot\nabla H_i(q_i)}{\lambda R_i^\beta}\right)
       -\mathcal F(\vartheta)
\]
have weighted $\ell^1_\rho$ norm and parameter Lipschitz seminorm at most
$C\lambda^{\min(2,4/s)}$. On $\Omega_i(4)$,
\begin{equation}\label{an:Hmacro}
 |\nabla H_i|\le C\lambda R_i,
 \qquad |D^2H_i|\le C\lambda.
\end{equation}
The differences after pullback by $T_i$ satisfy the same bounds multiplied
by $\|\vartheta-\vartheta'\|_X$.
\end{proposition}
\begin{proof}
By \eqref{an:projected-integral} and \eqref{an:delta-source},
\begin{equation}\label{an:Hidentity}
 \begin{split}
 H_i-B_i={}&Q_i-\Gamma*(\zeta_iQ_i^p)
 +\sum_{j\ne i}N_2(\Gamma*\delta f_j)\\
 &-T_2(\Gamma*\delta f_i)+\mathcal S f_{\rm out}.
 \end{split}
\end{equation}
The first term equals $\Gamma*((1-\zeta_i)Q_i^p)$ with the source taken
over the whole space, because $Q_i=\Gamma*Q_i^p$. Its normalized
$C^2$ jet on $\Omega_i(4)$ is $O((\mu_i/R_i)^4)$ after division by
$\lambda R_i^\beta$ with the appropriate derivative factors.

Set $M_i=\|\delta f_i\|_1$. The coarse bound
$|w_i|\le C\lambda R_i^2$ and the multiplier estimate give
\begin{equation}\label{an:mass-input}
 \frac{M_i}{\lambda R_i^{s+\beta}}
 \le C\left[R_i^{2-\beta}L_i^{-s}J_0(L_i)
       +\lambda^{p-1}R_i^{4+2p-\beta}+L_i^{-s}\right].
\end{equation}
The first term is bounded by
\[
 C\begin{cases}
 \lambda^2R_i^{s+\beta+2},&s<4,\\
 \lambda^2R_i^{6+\beta}(1+\log L_i),&s=4,\\
 \lambda^{8/s}R_i^{6-\beta+8\beta/s},&s>4.
 \end{cases}
\]
All shell exponents are strictly positive. By
\eqref{an:projected-lip} and the same mean value inequality, the actual
transported source difference, including its Jacobian, obeys the
right-hand side of \eqref{an:mass-input} multiplied by $\varepsilon$.

For sources distinct from $i$, their supports are separated from
$\Omega_i(4)$. Combine complete outer shells before applying the fourth-order
Taylor remainder, use the second kernel moment for inner shells, and use
macroscopic separation for the remaining centres on the target shell.
The resulting normalized $C^2$ jet matrix is precisely the bounded
matrix $H_{k\ell}$ from Theorem~\ref{an:countable}, applied to
$M_j/(\lambda R_j^{s+\beta})$. The own Taylor term obeys the same mass bound.
Its parameter difference uses the actual source difference just proved
and the kernel difference from Proposition~\ref{an:param}.

For the last term of \eqref{an:Hidentity}, the part
$(1-\sum\zeta_i)U^p$ is bounded by $C\lambda^p r^{2p}$ and vanishes
near every core. Splitting the potential into $|Y|<r/2$, comparable
scales, and $|Y|>2r$ gives
\[
 |D^j\mathcal S[(1-\sum\zeta_i)U^p](X)|
       \le C\lambda^p r^{4-j},\qquad 0\le j\le2.
\]
The absolute moments used here converge since $2p>0$.
The complete bubble-tail term has the same bound by
Lemma~\ref{anq:tails}. The actual differences satisfy these bounds
multiplied by $\varepsilon$; on the affine caps the kernel comparison
also applies to two derivatives.

Consequently the normalized jet of $H_i-B_i$ and its parameter
Lipschitz seminorm are bounded by
$C\lambda^{\min(2,8/s)}$; when $s=4$, include $1+|\log\lambda|$.
This is $O(\lambda^{\min(2,4/s)})$ for every $s\ge1$.
In particular, at $s=1$ the present mass estimate gives $O(\lambda^2)$;
the $O(\lambda^4)$ softening error alone does not determine the
effective-background error.
Add the softening error from Lemma~\ref{bg:bounds} to obtain the first
assertion. Since $R_i^{\beta-2}$ is bounded, the jet estimates and the
macroscopic background bounds also prove \eqref{an:Hmacro} and its difference.
\end{proof}

\subsection{Pohozaev identities and the reduced map}

For $y=X-q_i$, define
\begin{align}
 T_{jl}(w)&=w_j(\Delta w)_l-(\Delta w)w_{jl}
             +\tfrac12\delta_{jl}(\Delta w)^2,\notag\\
 D_l(w)&=y_jT_{jl}(w)+aw(\Delta w)_l-(a+1)(\Delta w)w_l.
 \label{an:stress}
\end{align}
Repeated indices are summed. Direct differentiation gives
$\partial_lT_{jl}=w_j\Delta^2w$ and
$\partial_lD_l=(y\cdot\nabla w+aw)\Delta^2w$.
Let $\mathcal B_i=\Omega_i(4)$, $e_i=q_i/|q_i|$, and put
\begin{align}
 P_{i,r}(w)&=\int_{\partial\mathcal B_i}
 \left((e_i)_jT_{jl}(w)\nu_l
       -\frac{e_i\cdot\nu}{p+1}w^{p+1}\right)\dd S,\notag\\
 P_{i,0}(w)&=\int_{\partial\mathcal B_i}
 \left(D_l(w)\nu_l-\frac{y\cdot\nu}{p+1}w^{p+1}\right)\dd S.
 \label{an:flux}
\end{align}
Since $a=n/(p+1)$, the nonlinear volume terms cancel and
\begin{equation}\label{an:flux-E}
 P_{i,r}(U)=\int\partial_{e_i}U\,\mathcal E_i,
 \qquad P_{i,0}(U)=\int(y\cdot\nabla U+aU)\mathcal E_i.
\end{equation}
All unspecified local integrals in this subsection are over $\mathcal B_i$.

For any smooth compactly supported real source $f$, including a signed
one, symmetry and the degree $-s$ of $\Gamma$ give
\begin{equation}\label{an:selfzero}
 \int f\,\partial_e(\Gamma*f)=0,
 \qquad \int f[(X-q)\cdot\nabla+a](\Gamma*f)=0.
\end{equation}
Indeed the first double integral changes sign on exchanging its variables;
symmetrizing the second gradient integral gives $-s\int f(\Gamma*f)/2$.
The kernel and its first derivatives are locally integrable, so both
exchanges are absolute. Apply these identities to $f_i$ and integrate the
remaining nonlinear terms by parts once more. Since $\zeta_i$ vanishes
near $\partial\mathcal B_i$ and $a=n/(p+1)$, the boundary terms cancel
exactly, and
\begin{align}
 P_{i,r}(U)&=\int\partial_{e_i}H_i f_i
 +\frac1{p+1}\int(\partial_{e_i}\zeta_i)U^{p+1},
 \label{an:mixed-r}\\
 P_{i,0}(U)&=\int(y\cdot\nabla H_i+aH_i)f_i
 +\frac1{p+1}\int(y\cdot\nabla\zeta_i)U^{p+1}.
 \label{an:mixed-0}
\end{align}
Indeed, differentiating $1-\zeta_i$ gives the positive cutoff term in
each formula; in the second, the remaining volume coefficient is
$a-n/(p+1)=0$. Denote these last two integrals, including the factor
$(p+1)^{-1}$, by $J_{i,r}$ and $J_{i,0}$. They are supported in
$\Omega_i(3)\setminus\Omega_i(2)$.

\begin{proposition}\label{an:reduction}
Define, using one centre from each shell,
\[
 \mathcal G_{\lambda,k}(\vartheta)=
 \left(\frac{P_{i,0}(U_\vartheta)}{a m_i\lambda R_i^\beta},
       \frac{|q_i|P_{i,r}(U_\vartheta)}{m_i\lambda R_i^\beta}\right).
\]
Then $\mathcal G_\lambda=\mathcal F+\mathscr R_\lambda$, where
\begin{equation}\label{an:reduced-error}
 \|\mathscr R_\lambda(0)\|_{\ell^1_\rho}
       +\operatorname{Lip}(\mathscr R_\lambda)
       \le C\lambda^{2/s}.
\end{equation}
\end{proposition}
\begin{proof}
Set $\varepsilon=\|\vartheta-\vartheta'\|_X$; a prime denotes the second
parameter. All constants below are independent of the shell and of
sufficiently small $\lambda$.

We first compare the cutoff terms on a fixed annulus. Write
$q_i^0=q_i(0)$, $\zeta_i^0=\zeta_i\circ T_i^\vartheta$, and
\[
 \mathcal A_i^0=B(q_i^0,3\delta R_i)\setminus B(q_i^0,2\delta R_i),
 \qquad \widetilde U_i^\vartheta=U_\vartheta\circ T_i^\vartheta.
\]
The cutoff $\zeta_i^0$ is parameter independent. The unit vector
$e_i=q_i/|q_i|$ is also fixed, since the centres move only radially.
The change of variables $X=T_i^\vartheta(Y)$ has Jacobian
$\sigma_i^n$ and gives
\begin{align}
 J_{i,r}^\vartheta
 &=\frac{\sigma_i^{n-1}}{p+1}
   \int_{\mathcal A_i^0}(\partial_{e_i}\zeta_i^0)
                  (\widetilde U_i^\vartheta)^{p+1}\dd Y,\notag\\
 J_{i,0}^\vartheta
 &=\frac{\sigma_i^n}{p+1}
   \int_{\mathcal A_i^0}((Y-q_i^0)\cdot\nabla\zeta_i^0)
                  (\widetilde U_i^\vartheta)^{p+1}\dd Y.
 \label{an:flux-rev-pullback}
\end{align}
Here $\sigma_i=e^{\eta_i/a}$, so $\sigma_i^{n-1}$ and $\sigma_i^n$ are
uniformly bounded, and their actual differences are at most $C\varepsilon$.
Lemma~\ref{bg:bounds} and \eqref{an:projected-bound}--\eqref{an:projected-lip}
give, on $\mathcal A_i^0$,
\[
 |\widetilde U_i^\vartheta|\le C\lambda R_i^2,
 \qquad |\widetilde U_i^\vartheta-\widetilde U_i^{\vartheta'}|
       \le C\varepsilon\lambda R_i^2.
\]
Apply the mean value inequality to the power $p+1$. Using
$|\nabla\zeta_i^0|\le C R_i^{-1}$ and $|\mathcal A_i^0|\le C R_i^n$,
we obtain
\begin{align}
 |J_{i,r}|&\le C\lambda^{p+1}R_i^{n+2p+1},&
 |J_{i,r}-J'_{i,r}|&\le C\varepsilon\lambda^{p+1}R_i^{n+2p+1},\notag\\
 |J_{i,0}|&\le C\lambda^{p+1}R_i^{n+2p+2},&
 |J_{i,0}-J'_{i,0}|&\le C\varepsilon\lambda^{p+1}R_i^{n+2p+2}.
 \label{an:flux-rev-transition}
\end{align}

The denominators also vary with the parameters. Put
$\widehat m_i=\int_{\mathbb R^n}\widehat f_i$.
Lemma~\ref{an:moments} and the relative difference of $\mu_i$ imply
\begin{equation}\label{an:flux-rev-mass}
 \begin{split}
 &\widehat m_i=I_p+o(1),\qquad
 |\widehat m_i-\widehat m_i'|\le\omega(\lambda)\varepsilon,\\
 &m_i=\mu_i^a\widehat m_i\asymp\lambda R_i^{s+\beta},\qquad
 |m_i-m_i'|\le C\varepsilon m_i,
 \end{split}
\end{equation}
where $\omega(\lambda)\to0$. In particular both microscopic masses are
bounded away from zero. Since $|q_i|\asymp R_i$ and
$\big||q_i|-|q_i'|\big|\le C\varepsilon R_i$, the two normalized cutoff
terms and their differences are bounded, respectively, by
\begin{equation}\label{an:flux-rev-normalized}
 C\lambda^{p-1}R_i^{2p+6-2\beta}
   =C\lambda^{8/s}R_i^{8-2\beta+16/s}
 \quad\hbox{and}\quad
 C\varepsilon\lambda^{8/s}R_i^{8-2\beta+16/s}.
\end{equation}
This includes the differences of $m_i^{-1}$ and of the radial prefactor
$|q_i|$.

For the source terms, use the common microscopic variable $y$ and define
\begin{align}
 K_{i,0}^\vartheta(y)
 &=H_i(q_i+\mu_i y)-H_i(q_i)
       +a^{-1}\mu_i y\cdot\nabla H_i(q_i+\mu_i y),\notag\\
 K_{i,r}^\vartheta(y)
 &=\partial_{e_i}H_i(q_i+\mu_i y)-\partial_{e_i}H_i(q_i),\notag\\
 \mathcal T_{i,\alpha}^\vartheta
 &=\frac{\int_{\mathbb R^n}K_{i,\alpha}^\vartheta(y)
                     \widehat f_i^\vartheta(y)\dd y}
          {\widehat m_i^\vartheta},\qquad \alpha\in\{0,r\}.
 \label{an:flux-rev-taylor}
\end{align}
The supports in this variable lie in the fixed ball $|y|\le3L_i$.
Proposition~\ref{an:background} gives, on that ball,
\begin{align*}
 |K_{i,0}^\vartheta(y)|&\le C\lambda\mu_iR_i|y|,&
 |K_{i,0}^\vartheta(y)-K_{i,0}^{\vartheta'}(y)|
       &\le C\varepsilon\lambda\mu_iR_i|y|,\\
 |K_{i,r}^\vartheta(y)|&\le C\lambda\mu_i|y|,&
 |K_{i,r}^\vartheta(y)-K_{i,r}^{\vartheta'}(y)|
       &\le C\varepsilon\lambda\mu_i|y|.
\end{align*}
Indeed, $T_i^\vartheta(q_i^0+\mu_i(0)y)=q_i+\mu_i y$.
For the scale bracket the factors $\sigma_i$ in the displacement and
$\sigma_i^{-1}$ in the pulled-back gradient cancel. For the radial
bracket there remains the bounded factor $\sigma_i^{-1}$, whose difference
is $O(\varepsilon)$. The common-coordinate $C^2$ estimates therefore
give the stated differences, with either of the comparable scales
$\mu_i,\mu_i'$ on the right.

For completeness, the actual difference of each quotient is
\begin{equation}\label{an:flux-rev-quotient}
 \begin{split}
 \mathcal T_{i,\alpha}^\vartheta-\mathcal T_{i,\alpha}^{\vartheta'}
 ={}&\frac{\int(K_{i,\alpha}^\vartheta-K_{i,\alpha}^{\vartheta'})
                                      \widehat f_i^\vartheta}
               {\widehat m_i^\vartheta}\\
 &+\frac{\int K_{i,\alpha}^{\vartheta'}
                     (\widehat f_i^\vartheta-\widehat f_i^{\vartheta'})}
               {\widehat m_i^\vartheta}\\
 &+\left(\frac1{\widehat m_i^\vartheta}
            -\frac1{\widehat m_i^{\vartheta'}}\right)
                 \int K_{i,\alpha}^{\vartheta'}\widehat f_i^{\vartheta'}.
 \end{split}
\end{equation}
The first absolute source moment controls the first term. Equation
\eqref{an:micro-source-lip} and \eqref{an:flux-rev-mass} control the other
two. Each retains the factor $\mu_i$ from $K_{i,\alpha}$; hence the
modulus $\omega(\lambda)=o(1)$ need not have the rate $\lambda^{2/s}$.
It follows that $\mathcal T_{i,0}$ and $\mathcal T_{i,r}$ have bounds
$C\lambda\mu_iR_i$ and $C\lambda\mu_i$, with the same bounds multiplied
by $\varepsilon$ for their actual differences.

The exact identities now read
\[
 \frac{P_{i,0}}{a m_i}=H_i(q_i)+\mathcal T_{i,0}+\frac{J_{i,0}}{a m_i},
 \qquad
 \frac{P_{i,r}}{m_i}=\partial_{e_i}H_i(q_i)
                     +\mathcal T_{i,r}+\frac{J_{i,r}}{m_i}.
\]
Thus both normalized Taylor errors and their Lipschitz constants are
bounded by
\[
 C\mu_iR_i^{1-\beta}
       \le C\lambda^{2/s}R_i^{3-\beta+2\beta/s}.
\]
The exponents $\gamma_1=3-\beta+2\beta/s$ and
$\gamma_2=8-2\beta+16/s$ are positive by \eqref{an:parameter-range}.
For $\gamma_2$ one only needs $\beta<4$.
In the choice of the weight take $\rho>1$ sufficiently close to one that
$\rho t^{\gamma_j}<1$ for $j=1,2$. Since $R_k=R_0t^k$,
\begin{equation}\label{an:flux-rev-sum}
 \sum_{k\ge0}\rho^kR_k^{\gamma_j}
       =\frac{R_0^{\gamma_j}}{1-\rho t^{\gamma_j}}<\infty.
\end{equation}
Weighted summation and Proposition~\ref{an:background} consequently give
\[
 \|\mathscr R_\lambda(0)\|_{\ell^1_\rho}
       +\operatorname{Lip}(\mathscr R_\lambda)
 \le C(\lambda^{\min(2,4/s)}+\lambda^{2/s}+\lambda^{8/s})
 \le C\lambda^{2/s}.
\]
The last inequality also holds at $s=1$: the three terms have orders
$\lambda^2$, $\lambda^2$ and $\lambda^8$. The shell exponents in the
two flux errors are then $3+\beta$ and $24-2\beta$, both positive.
The factor $|q_i|$ in the radial row identifies its principal term
exactly with $q_i\cdot\nabla H_i(q_i)$, completing the comparison with
$\mathcal F$.
\end{proof}

\section{Solving the parameters and removing the multipliers}

\begin{theorem}\label{an:exact}
For sufficiently small $\lambda>0$, there is a parameter
$\vartheta_\lambda\in X$, with $\|\vartheta_\lambda\|_X\le C\lambda^{2/s}$,
such that every multiplier in Proposition~\ref{an:projection} vanishes.
Consequently $U=W_{\vartheta_\lambda}+\phi_{\vartheta_\lambda}$ is positive and
smooth in $B_1\setminus\{0\}$ and satisfies $\Delta^2U=U^p$ there.
\end{theorem}
\begin{proof}
Let $A=\mathcal F'(0)$ and let $\mathcal R$ be its bounded right inverse
from Proposition~\ref{sym:rightinverse}. Its range is closed, since
$\mathcal R A$ is a bounded projection. Set
$\mathcal Q(\vartheta)=\mathcal F(\vartheta)-A\vartheta$ and solve in this range
\[
 \vartheta=-\mathcal R[\mathcal Q(\vartheta)+\mathscr R_\lambda(\vartheta)].
\]
On the ball $\|\vartheta\|_X\le K\lambda^{2/s}$, choose
$K>2C\|\mathcal R\|$, and then $\lambda$ small.
Proposition~\ref{sym:regularity} and \eqref{an:reduced-error} show that this
map preserves the ball and has Lipschitz constant at most
$C\|\mathcal R\|(K+1)\lambda^{2/s}<1/2$.
Its fixed point satisfies $\mathcal G_\lambda=0$, because
$A\mathcal R=I$. The finitely many varied exterior sources remain outside
$B_1$ and have positive masses. No equation for an exterior core is imposed.

It remains to prove that zero flux eliminates the multipliers.
Write $w_i=U-Q_i$ and, in microscopic coordinates, set
$\psi_{\alpha,L}=\chi_LVZ_\alpha$.
The functions $VZ_0$ and $VZ_j$ decay respectively as
$(1+|y|)^{-s-8}$ and $(1+|y|)^{-s-9}$; their derivatives gain one power.
Hence, uniformly for $L\ge2$,
\begin{equation}\label{an:test-moments}
 \int(1+|y|^2)
 (|\psi_{\alpha,L}|+|y\cdot\nabla\psi_{\alpha,L}|
                         +|\nabla\psi_{\alpha,L}|)\dd y\le C.
\end{equation}
The cutoff derivative terms are bounded on $L/2<|y|<L$ by constants
times $|VZ_\alpha|$ and $L^{-1}|VZ_\alpha|$; with an extra $|y|^j$,
their integrals are $O(L^{j-4})$ and $O(L^{j-5})$, respectively, for $j\le2$.

For $b_i=C\lambda R_i^2$, we have $|w_i|\le b_i$ on the multiplier support.
Integration by parts, with no boundary term, gives
\begin{align*}
 \int\chi_iQ_i^{p-1}Z_{i,\alpha}((X-q_i)\cdot\nabla+a)w_i
 &=\mu_i^s\int[(a-n)\psi_{\alpha,L}
             -y\cdot\nabla\psi_{\alpha,L}]w_i(q_i+\mu_i y)\dd y,\\
 \int\chi_iQ_i^{p-1}Z_{i,\alpha}\partial_jw_i
 &=-\mu_i^{s-1}\int\partial_j\psi_{\alpha,L}
                               w_i(q_i+\mu_i y)\dd y.
\end{align*}
The respective errors are $O(\mu_i^sb_i)$ and $O(\mu_i^{s-1}b_i)$.
For the principal bubble,
$( (X-q_i)\cdot\nabla+a)Q_i=\mu_i^{-a}Z_{i,0}$ and
$\partial_jQ_i=\mu_i^{-a-1}Z_{i,j}$.

Keep the original, equally normalized unknowns $C_{i,\alpha}$.
Divide the zero dilation-flux equation by $\mu_i^a$ and each zero
translation-flux equation by $\mu_i^{a-1}$. The resulting coefficient
matrix is
\[
 G_{L_i}+O(\mu_i^ab_i),\qquad
 (G_L)_{\alpha\gamma}=\int\chi_LVZ_\alpha Z_\gamma.
\]
Radiality of $\chi_L$ makes the mixed entries zero, and
$G_L\to g_nI$. Moreover,
$\mu_i^ab_i\le C\lambda^2R_i^{s+\beta+2}\to0$ uniformly in $i$.
Thus these matrices are uniformly invertible. In the invariant class only
the dilation and radial translation equations remain, so the same conclusion
applies to that two-by-two block. Equation~\eqref{an:flux-E} gives
$C_{i,\alpha}=0$ for every inner centre. Proposition~\ref{an:projection}
now proves the theorem.
\end{proof}

\section{The singularity and the spherical lower envelope}\label{end:section}

\begin{proof}[Proof of Theorem~\ref{intro:main}]
Apply Theorem~\ref{an:exact} with the geometry and approximate solution constructed
above. It gives parameters $\vartheta_\lambda$ with
$\|\vartheta_\lambda\|_X\le C\lambda^{2/s}$ and a positive exact solution
\[
U=W_{\vartheta_\lambda}+\phi,\qquad
\Delta^2U=U^p\quad\text{in }B_1\setminus\{0\}.
\]
Every nonzero center is a smooth point; elliptic regularity and
positivity give $U\in C^\infty(B_1\setminus\{0\})$.
Moreover $\|\phi\|_{\mathcal Y}\le C\lambda^{1+2/s}$.
At the centers,
\begin{equation}\label{end:peaks}
U(q_{i,\nu})=A_n\mu_i^{-a}+O(\lambda R_i^2)
           \ge c\lambda^{-1}R_i^{-s-\beta}\longrightarrow\infty.
\end{equation}

The normal subspace is separated from every cap. The point-field
bounds and the soft-error estimate give, for $r=|z|$,
\begin{equation}\label{end:axis-W}
W_{\vartheta_\lambda}(0,z)
 =\lambda\kappa|z|^2+O(\lambda r^\beta)
                         +O(\lambda^{1+4/s}r^4).
\end{equation}
Since $|\phi(0,z)|\le C\|\phi\|_{\mathcal Y}r^\beta$, we have
$U(0,z)\to0$. Positivity and \eqref{end:peaks} prove
\eqref{intro:oscillation}.

We must also control two derivatives of the correction; its
$\mathcal Y$ norm alone does not give this. In the source decomposition
of Sect.~\ref{an:effective-section}, let
$M_i=\|\delta f_i\|_{L^1}$. The mass estimates there imply
$M_i\le C\lambda R_i^{s+\beta}$.
Here and below the finite number of vertices on each shell is
absorbed in $C$. At $(0,z)$ the distance to the $i$th source is at least
$c\max(r,R_i)$. For outer shells first combine the complete even
orbit, and then differentiate the renormalized kernel. By
\eqref{bg:kernel},
\begin{align}
|D^2\phi(0,z)|
&\le C\lambda\left[
 r^2\sum_{R_i\gtrsim r}R_i^{\beta-4}
      +\sum_{R_i\lesssim r}R_i^{\beta-2}\right]
       +C\lambda^p r^2\notag\\
&\le C\lambda r^{\beta-2}+C\lambda^p r^2\longrightarrow0.
\label{end:axis-phi}
\end{align}
Comparable shells are covered by the same bound, since the observation
point remains separated from their caps. The last term comes from the
outer source $f_{\rm out}$, whose renormalized potential has second
derivatives $O(\lambda^p r^2)$.

The same shell estimates and \eqref{bg:soft-sums} give
\[
D^2(W_{\vartheta_\lambda}-\lambda\kappa|z|^2)(0,z)
 =O(\lambda r^{\beta-2})+O(\lambda^{1+4/s}r^2).
\]
Consequently,
\begin{equation}\label{end:limit}
\lim_{z\to0,\ z\ne0}(-\Delta U)(0,z)
      =-2(n-2)\lambda\kappa<0.
\end{equation}
Set $u=U$ and $c_*=2(n-2)\lambda\kappa$.
This proves \eqref{intro:negative}. The unbounded sequence
\eqref{end:peaks} shows that the origin is nonremovable.
There are no other singularities. Finally, a radial function would
take the same value at $q_{i,\nu}$ and at a normal point of radius
$|q_{i,\nu}|$, contradicting \eqref{end:peaks} and
\eqref{end:axis-W}. The theorem follows.
\end{proof}

\begin{proof}[Proof of Corollary~\ref{intro:envelope}]
Fix $\lambda>0$ sufficiently small. By Lemma~\ref{bg:bounds} and the
smallness of $\|\phi\|_{\mathcal Y}/\lambda$,
\[
U(X)\ge c\lambda|X|^\beta.
\]
Choose any unit vector $v_*\in S^1\setminus V$.
Its angular distance from the finite vertex set is positive.
For every $r>0$ and every source center,
\[
|(rv_*,0)-q_{i,\nu}|\ge c\max(r,R_i).
\]
For comparable radii this follows from the fixed angular gap; for
noncomparable radii it follows from radial separation.
The renormalized shell estimates thus give
\[
|F_{\vartheta_\lambda}(rv_*,0)|\le Cr^\beta.
\]
The normal polynomial
vanishes at these points, and the soft error is
$O(\lambda^{1+4/s}r^4)$. Hence
\[
U(rv_*,0)\le C\lambda r^\beta.
\]
This proves \eqref{intro:minimum}, with $\lambda$ absorbed into the
constants.

Take $r_i=|q_{i,\nu}|\asymp R_i$. Equation \eqref{end:peaks} and the
upper bound for the spherical minimum give
\[
M_U(r_i)\ge c\lambda^{-1}r_i^{-s-\beta},\qquad
\frac{M_U(r_i)}{m_U(r_i)}
       \ge c\lambda^{-2}r_i^{-s-2\beta}.
\]
Again absorb the fixed $\lambda$ into the constants. If
$U(X)\ge b|X|^\gamma$ for $\gamma<\beta$, then
$b\le Cr^{\beta-\gamma}$ for every small $r$, a contradiction.
Finally, $s+\beta>s/2$, so the peak estimate excludes the critical
upper bound in the statement.
\end{proof}

\medskip

\noindent{\bf Acknowledgments.} This work was completed while H. Yang was visiting The Chinese University of Hong Kong. He would like to thank the Department of Mathematics and the Institute of Mathematical Sciences for their hospitality. 
The research of J. Wei was supported by National R\&D Program of China (Grant No. 2022YFA1005602), and Hong Kong General Research Fund (No. 14303125) ``On Fujita equation in the critical or supercritical regime”. The research of H. Yang is supported by NSFC 12301140. The authors used AI-assisted tools, principally ChatGPT/Codex, for brainstorming and manuscript editing. AI output was not treated as a mathematical authority. The authors independently worked and verified all theorem statements, proofs, calculations, constants, citations, and final prose, and take full responsibility for the contents of the paper.

\noindent Juncheng Wei 

\noindent Department of Mathematics \\ 
The Chinese University of Hong Kong, Hong Kong \\[1mm]
Email: \textsf{wei@math.cuhk.edu.hk}

\medskip 

\noindent Hui Yang

\noindent School of Mathematical Sciences \\
Shanghai Jiao Tong University, Shanghai 200240, China \\[1mm]
Email: \textsf{hui-yang@sjtu.edu.cn}

\medskip 

\noindent Ronghao Yang

\noindent School of Mathematical Sciences \\
Shanghai Jiao Tong University, Shanghai 200240, China \\[1mm]
Email: \textsf{ylqfyrh5218@sjtu.edu.cn}

\end{document}